\documentclass{article}
\usepackage{amssymb}
\usepackage{epsfig}
\usepackage{amsmath}
\usepackage{amsfonts}

\usepackage{color}
\usepackage{hyperref}

\newcommand{\N}{{{\Bbb N}}}

\newtheorem{theorem}{\sc Theorem}[section]
\newtheorem{proposition}{\sc Proposition}[section]
\newtheorem{lemma}{\sc Lemma}[section]
\newtheorem{definition}{\sc Definition}[section]
\newtheorem{remark}{\sc Remark}[section]
\newtheorem{corollary}{\sc Corollary}[section]
\newtheorem{example}{\sc Example}[section]

\title{Coupled fixed points of multivalued operators and first--order ODEs with state--dependent deviating arguments\footnote{Partially
supported by FEDER and
Ministerio de Educaci\'on y Ciencia, Spain,
project MTM2010--15314.}}  

\author{Rub\'en Figueroa and Rodrigo L\'opez Pouso}

\begin{document}
\maketitle

\begin{center}
 Departamento de
An\'alise Matem\'atica\\
Facultade de Matem\'aticas\\Universidade de Santiago de Compostela, Campus Sur \\
15782 Santiago de
Compostela, Spain\\
  {\bf e--mail:} ruben.figueroa@usc.es, 
   rodrigo.lopez@usc.es
\end{center}

\def\tanh#1{\,{\normalsize tanh}{\,#1}\,}
\def\simbolo#1#2{#1\dotfill{#2}}
\def\qed{\hbox to 0pt{}\hfill$\rlap{$\sqcap$}\sqcup$\medbreak}
\def\theequation{\arabic{section}.\arabic{equation}}
\def\thesection {\arabic{section}}

\begin{abstract}
We establish a coupled fixed points theorem for a meaningful class of mixed monotone multivalued operators and then we use it to derive some results on existence of quasisolutions and solutions to first--order functional differential equations with state--dependent deviating arguments. Our results are very general and can be applied to functional equations featuring discontinuities with respect to all of their arguments, but we emphasize that they are new even for differential equations with continuously state--dependent delays.
\end{abstract}

\noindent
{\bf Primary classification number:} 34K05.
 
\section{Introduction}

We present an abstract result on coupled fixed points for multivalued operators and we use it to prove the existence of absolutely continuous solutions to the problem
\begin{equation}\label{p1}
\left\{
\begin{array}{ll}
x'(t)=f(t,x(t),x(\tau(t,x(t),x)),x)   \quad \mbox{for almost all $t \in I_+=[t_0,t_0+L],$} \\
\\
x(t)=\Lambda(x)+k(t) \quad \mbox{for all $t \in I_-=[t_0-r,t_0],$}
\end{array}
\right.
\end{equation}
where $t_0 \in \mathbb R$, $L>0$, $r \ge 0$, $k$ is a given continuous function when $r>0$ or a constant in case $r=0$.

Let ${\cal C}(I_{\pm})$ denote the space of continuous functions on $I_{\pm}=[t_0-r,t_0+L]$, and let $AC(I_+)$ denote the space of absolutely continuous functions on $I_+$. The functions 
$$f:Dom(f) \subset I_+ \times \mathbb R^2 \times {\cal C}(I_{\pm})  \longrightarrow \mathbb R,$$
$$ \tau:Dom(\tau) \subset I_+ \times \mathbb R \times {\cal C}(I_{\pm})  \longrightarrow I_{\pm},$$
 and $\Lambda:Dom(\Lambda)\subset {\cal C}(I_{\pm})
\to \mathbb R$ need not be continuous with respect to some of their arguments, as we will specify later, but $f$ will be, roughly speaking, nonincreasing with respect to its third variable and nondecreasing with respect to its fourth one. The deviating function $\tau$ can assume either past or future values at each point $(t,x,\phi) \in Dom(\tau)$.
 \bigbreak

Differential equations with deviating arguments are being intensively studied, and we can quote recent papers as \cite{dyk}, \cite{gao}, \cite{jan2}. This paper considers the general situation where deviations depend on the unknown, which is interesting from the viewpoints of both theory and applications. Indeed, the particular case of differential equations with state--dependent delays is receiving a lot of attention, see for instance \cite{nto}, \cite{her1}, \cite{wal1}, \cite{wal2}, \cite{zen}. We also refer readers to the survey by Hartung {\it et al.} \cite{har} and references therein, for mathematical models which use state--dependent delays, and for basic and qualitative theory on this type of problems.
\bigbreak
Our starting point is the work by Jankowski in \cite{jan}, who studied the problem
\begin{equation}\label{janp}
x'(t)=f(t,x(\tau(t,x(t)))) \mbox{ for all $t \in I_+$,} \quad  x(t_0)=\lambda x(t_0+L) + k \quad (\lambda, \, k \in \mathbb R),
\end{equation}
for continuous $f$ and $\tau$ and monotonicity with respect to spatial variables. In that paper, the author uses a monotone method in presence of lower and upper solutions in order to obtain the existence of extremal quasisolutions for problem (\ref{janp}) and then he applies a maximum principle to guarantee that, under certain Lipschitz conditions, those extremal quasisolutions are the same function, thus proving that problem (\ref{janp}) has a unique solution between the given lower and upper solutions. The author also assumes in \cite{jan} that $f \ge 0$ on the compact subset of $I_+ \times \mathbb R$ delimited by the graphs of the lower and the upper solutions.  \\

It is the aim of the present work to state and prove a new result on coupled fixed points for multivalued operators and then apply it to (\ref{p1}), showing in this way that the existence results in \cite{jan} hold valid in more general settings, namely,
\begin{enumerate}
\item In the case $\tau(t,x(t),x)=\tau(t,x(t))$ the nonlinearity $f$ has to be neither nonnegative nor monotone with respect to $x(\tau)$ between assumed lower and upper solutions;
\item Delay differential equations are included in the framework of (\ref{p1});
\item Neither $f$ nor $\tau$ need be continuous in any of their arguments for the existence of quasisolutions, and they can be discontinuous with respect to the $t$ and the $x(t)$ variables for the existence of solutions;
\item The right--hand side in the differential equation features a nondecreasing functional dependence in its fourth variable which makes nonincreasingness with respect to the third variable be less stringent, as we show in Section 3.4. Moreover, it provides us with a unified framework for the study of many other types of functional differential equations (integro--differential equations, equations with maxima, and so on);
\item The deviating function $\tau$ depends at each moment $t$ on the global behavior of the solution, and not merely on the value the solution assumes at $t$. Moreover, $\tau$ is not necessarily monotone with respect to $x(t)$;
    \item The ``boundary" conditions are allowed to be nonlinear and they also involve the global behavior of the unknown, and not only the values it assumes on $t_0$ and on $t_0+L$.
\end{enumerate}

Our approach uses a general fixed point theorem for nondecreasing ope\-ra\-tors in ordered spaces which can be looked up in \cite{hela}, combined with maximum principles and some other ideas originally introduced in \cite{jan}. Another basic result in this research is a recent existence theorem for first--order discontinuous differential equations proven by the authors in \cite{figp}, which we include in Section $3.1$ for the sake of self--containedness.
\\

We organize this paper as follows: in Section 2 we state and prove a new abstract result on existence of coupled fixed points for multivalued operators. Section 3 contains two new theorems on existence of coupled quasisolutions for problem (\ref{p1}) which, in turn, yield two corresponding new theorems on existence of unique solutions between given lower and upper solutions. We also include in Section 3 a complete discussion about the assumptions introduced in it. In Section 4 we present some sufficient conditions for existence of lower and upper solutions, thus facilitating the applicability of the results in Section 3. Finally, Section 5 considers  the particular case of problem (\ref{p1}) when $\tau(t,x(t),x)=\tau(t,x(t))$, and we show that, unlike in the general case considered in Section 3, no restriction on the sign of $f$ is needed. We illustrate the applicability of our results with some examples.

\section{Coupled fixed points for multivalued operators}
The following concepts and results are taken from \cite{hela}, although the notion of coupled (quasi)fixed points goes back to Moore \cite{moore} and the state of the art on coupled fixed points for mixed monotone operators owes lots of results and applications to Guo and Lakshmikantham, starting in \cite{gl}, see also \cite{gl2}.

\bigbreak

A metric space $X$ equipped with a partial ordering is an ordered metric space if the intervals $[x)=\{y \in X \, :\, x \le y\}$ and $(x]=\{y \in X \, :\, y \le x\}$ are closed for every $x \in X$.

Let $P$ be a subset of an ordered metric space; an operator $A:P \times P \longrightarrow P$ is mixed monotone if $A(\cdot,z)$ is nondecreasing and $A(z,\cdot)$ is nonincreasing for each $z \in P$. We say that $A$ satisfies the mixed monotone convergence property (m.m.c.p.) if $(A(v_j,w_j))_{j=1}^{\infty}$ converges whenever $(v_j)_{j=1}^{\infty}$ and $(w_j)_{j=1}^{\infty}$ are sequences in $P$, one being nondecreasing and the other nonincreasing.

\medbreak

Next we introduce the concept of extremal coupled fixed points for multivalued operators.
\begin{definition}
Let $P$ be a subset of an ordered metric space. We say that $v,w \in P$ are coupled fixed points of a multivalued operator ${\cal A}:P\times P \longrightarrow 2^{P}\backslash \emptyset$ if they fulfill

\medbreak

\begin{enumerate}
\item[$(CFP)$]{ $v \in {\cal A}(v,w)$ and $w \in {\cal A}(w,v)$}.
\end{enumerate}
We say that $v_*,w^*$ are the extremal coupled fixed points of ${\cal A}$ in $P$ if $v_*,w^*$ satisfy $(CFP)$ and

\medbreak

\begin{enumerate}
\item[$(ECFP)$]{ if $v,w \in P$ satisfies $(CFP)$, then $v_* \leq v $ and $w \leq w^*$. }
\end{enumerate}
\end{definition}
Notice that if $v_*,w^*$ are the extremal coupled fixed points of ${\cal A}$ in $P$ then $(ECFP)$ implies that $v_* \leq w^*$ because $w^*,v_*$ satisfy $(CFP)$.
\medbreak
For the convenience of the reader, we state a lemma on nondecreasing operators in ordered metric spaces which we use in the proof of our new result on coupled fixed points for multivalued operators.

\begin{lemma}\label{hela2}{\bf \cite[Theorem 1.2.2]{hela}} Let $Y$ be a subset of an ordered me\-tric space $X$, $[a,b]$ a nonempty closed interval in $Y$, and $G:[a,b] \longrightarrow [a,b]$ a nondecreasing operator.

If $(Gx_m)_{m=1}^{\infty}$ converges whenever $(x_m)_{m=1}^{\infty}$ is a monotone sequence  then $G$ has the extremal fixed points in $[a,b]$, i.e., a least fixed point, $x_*$, and a greatest one, $x^*$. Moreover
\begin{equation*}
x_*=\min\{y \, : \, Gy \leq y \} \quad \mbox{and} \quad  x^*=\max\{y \, : \, y \leq Gy \}.
\end{equation*}
\end{lemma}

As a consequence of Lemma \ref{hela2} we have the following result on existence of extremal coupled fixed points for multivalued mixed monotone operators. It can be regarded as a multivalued version of \cite[Theorem 1.2.4]{hela}.

\begin{theorem}\label{multi} Let $Y$ be a subset of an ordered me\-tric space $X$, $[\alpha,\beta]$ a nonempty closed interval in $Y$, and ${\cal A}:[\alpha,\beta] \times [\alpha,\beta] \longrightarrow 2^{[\alpha,\beta]} \backslash \emptyset$ a multivalued operator.

If for all $v,w \in [\alpha,\beta]$ there exist $A_-(v,w)=\min {\cal A}(v,w)$ and $A_+(v,w)=\max {\cal A}(v,w)$, and the (single--valued) operators $A_{\pm}$ are mixed monotone and satisfy the m.m.c.p., then ${\cal A}$ has the extremal coupled fixed points in $[\alpha,\beta]$.
\end{theorem}

\noindent {\bf Proof.}  In the product metric space $X \times X$ we consider the order $\preceq$ defined as follows: for two pairs $(v,w),(\overline{v},\overline{w}) \in X \times X$, we write $(v,w) \preceq (\overline{v},\overline{w})$ if $v \leq  \overline{v}$ and $w \geq \overline{w}$.

Let us denote $a=(\alpha,\beta)$ and $b=(\beta,\alpha)$, and consider the interval
$$[a,b]=\{(v,w) \in Y \times Y \, : \, a \preceq (v,w) \preceq b\},$$ and the operator $G: [a,b] \longrightarrow [a,b]$ which maps each pair $(v,w)$ into $$G(v,w)=(A_-(v,w), A_+(w,v)).$$

The assumptions readily imply that the conditions of Lemma \ref{hela2} are fulfilled, and therefore $G$ has the least fixed point $(v_*,w^*)$ in $[a,b]$, which, moreover, satisfies
\begin{equation}
\label{minmax}
(v_*,w^*)=\min_{\preceq}\{(v,w) \, : \, G(v,w) \preceq (v,w) \}.
\end{equation}

The definition of $G$ guarantees that $v_*$ and $w^*$ satisfy $(CFP)$. Moreover, if $(v,w) \in [a,b]$ satisfies $(CFP)$ then
$$A_-(v,w) \leq v \ \mbox{ and } A_+(w,v) \geq w,$$
so $G(v,w) \preceq (v,w)$ and then (\ref{minmax}) yields $(v_*,w^*) \preceq (v,w)$, or, equivalently, $v_* \leq v$ and $w^* \geq w$. Therefore the pair $(v_*,w^*)$ satisfies $(ECFP)$ with $P=[\alpha,\beta]$, and then $v_*,w^*$ are the extremal coupled fixed points of ${\cal A}$ in $[\alpha,\beta]$.
\qed

\section{Existence of quasisolutions and unique solutions}

Now we are ready to approach our original questions about the existence of quasisolutions and solutions for problem (\ref{p1}). We open this section with some preliminaries, and we close it with a discussion about the assumptions introduced in it.
\subsection{Preliminary results on discontinuous equations}
We recall some concepts and results from viability theory. The reader is referred to \cite{aucell}, \cite{cp} for more details.
\medbreak
For a
given set $A \subset \mathbb R^n$, $n \in \N$, the
Bouligand's contingent cone at $x
\in A$ is defined
as
$$T_A(x)=\bigcap_{\varepsilon >0} \bigcap_{\alpha>0}
\bigcup_{0<h<\alpha}
\left( \frac{1}{h}(A-x)+\varepsilon \, B \right),$$
where $B=(-1,1)$.

For an interval $I \subset \mathbb R$ and a multivalued map
$K:I \to 2^{\mathbb R}\backslash \emptyset $ the graph of $K$ is the
set
$graph(K)=\left\{ (t,x) \in I \times \mathbb R \, : \, x \in K(t)
\right\}$. The contingent derivative of $K$ at a point $(t,x)
\in graph(K)$ is defined as
the mapping $DK(t,x):\mathbb R \to 2^{\mathbb R}$ whose graph is the contingent
cone $T_{graph(K)}(t,x)$,
i.e.,
$$v_0 \in DK(t,x)(t_0) \Leftrightarrow
(t_0,v_0) \in T_{graph(K)}(t,x).$$

In case $K$ is single--valued its contingent derivative can be easily computed in the following significant
case:
\begin{proposition}
\label{dercont} \cite[Lemma 1.2]{cp} Let $I \subset \mathbb R$ be an
interval and
$\gamma:I \to \mathbb R$ a given function.

The mapping $K(t)=\{ \gamma(t)
\}$, $t
\in I$, satisfies that
$DK(t,\gamma(t))(1)$ lies between $D_+ \gamma(t)$ and
$D^+ \gamma(t)$ for
all $t\in I$, where $D_+ \gamma$ and $D^+
\gamma$ denote the lower--right and
the upper--right Dini
derivatives, respectively.

In particular, if $\gamma$
is right--differentiable at some $t\in
I$ then we have
$$DK(t,\gamma(t))(1)=\{ \gamma'_+(t)
\}.$$
\end{proposition}

Finally, we state an existence result for initial value problems with discontinuous ODEs which we extend in next sections in order to cover (\ref{p1}).

\begin{lemma}\cite[Theorem 3.1]{figp}
\label{lemdis}
Suppose that for the scalar initial value pro\-blem
\begin{equation}
\label{auxivp}
z'(t)=g(t,z(t)) \quad \mbox{for a.a. $t \in I_+$, \quad $z(t_0)=z_0 \in \mathbb R$,}
\end{equation}
the following conditions hold:
\medbreak
\begin{enumerate}
\item{(Lower and upper solutions) There exist  $\alpha, \beta \in AC(I_+)$ such
that $\alpha \leq \beta$ on $I_+$, the mappings $t\in I_+ \to g(t,
\alpha(t))$ and $t \in I_+ \to g(t,
\beta(t))$ are measurable,
and
$$
\alpha'(t) \leq g(t,\alpha(t)) \quad \mbox{for a.a. $t \in I_+$}, \quad \alpha(t_0) \leq
z_0,
$$
$$
\beta'(t) \geq g(t,\beta(t)) \quad \mbox{for a.a. $t \in I_+$}, \quad \beta(t_0) \geq
z_0;
$$}
\item{($L^1$ bound) There exists $\psi
\in L^1(I)$ such that
for a.a. $t \in I_+$ and for all $x \in [\alpha(t),\beta(t)]$ we
have $
|g(t,x)| \leq \psi(t);$}
\medbreak
\item{(Measurability) For all $x \in [\min_{t \in I_+}\alpha(t),\max_{t\in I_+}\beta(t)]$, the mapping
$$t \in \{s \in I_+ \, : \, \alpha(s) \le x \le \beta(s)\} \longmapsto g(t, x)$$ is
measurable;}
\medbreak
\item{(Admissible discontinuity sets) For a.a. $t \in I_+$, $g(t,\cdot)$ is
continuous on $[\alpha(t),\beta(t)]
\backslash K(t)$, where $K(t)= \cup_{n=1}^{\infty}K_n(t),$
and for
each $n \in \mathbb N$ and each $x \in K_n(t)$ we have
\begin{equation}
\label{unvres}
\bigcap_{\varepsilon >0}
\overline{{\rm co}} \, \, g(t, x + \varepsilon B)
\bigcap DK_n(t,x)(1) \subset
\{g(t,x)\},
\end{equation}
where $B=(-1,1)$ and $\overline{{\rm co}}$ denotes the closed convex hull.}
\end{enumerate}
\medbreak
Then the set of absolutely continuous solutions of problem (\ref{auxivp}) between $\alpha$ and $\beta$ on $I_+$
is nonempty and it has pointwise maximum, $z^*$, and mi\-ni\-mum,  $z_*$, which are the extremal solutions of (\ref{auxivp}) between $\alpha$ and $\beta$, and which satisfy
\begin{equation}\label{greatest}
z^*=\max \{v \in AC(I_+) \, : \, v'(s) \leq g(s,v(s)) \
\mbox{a.e., $v(t_0) \leq z_0$, $\alpha \leq v \leq \beta$}  \},
\end{equation}
\begin{equation}\label{least}
z_*=\min \{v \in AC(I_+)
\, : \, v'(s) \geq g(s,v(s)) \ \mbox{a.e., $v(t_0) \geq z_0$, $\alpha \leq v \leq \beta$}
\}.
\end{equation}

\end{lemma}

\begin{remark}
The assumptions in \cite[Theorem 3.1]{figp} concern all values of $x \in \mathbb R$, instead $x \in [\min_{t \in I_+}\alpha(t),\max_{t\in I_+}\beta(t)]$ or $x \in [\alpha(t),\beta(t)]$ in the relevant places. A revision of its proof shows that the result holds valid as stated in Lemma \ref{lemdis}.
\end{remark}

\begin{remark}
In the conditions of Lemma \ref{lemdis}, for $(t,x) \in I_+ \times \mathbb R$ we have
$$
\bigcap_{\varepsilon >0}\overline{{\rm co}} \, \, g(t, x + \varepsilon B)= \left[\min\left\{\liminf_{y \to x}g(t,y),g(t,x) \right\},
\max\left\{\limsup_{y \to x}g(t,y),g(t,x) \right\} \right].
$$
In most applications of Lemma \ref{lemdis} the multivalued maps $K_n$ can be chosen to be single--valued and then we can use Proposition \ref{dercont} to compute the contingent derivative when checking condition (\ref{unvres}). This is an advantage provided by the decomposition of the discontinuity sets $K(t)$ into countable unions. Another advantage of that decomposition has to do with the fact that we hope that intersections in (\ref{unvres}) be small when we use Lemma \ref{lemdis}, and clearly $DK_n(t,x)(1) \subset DK(t,x)(1)$.
\end{remark}

The next two sections contain our main contributions for problem (\ref{p1}). Note that we will often shorten notation by writing $\tau_{t,x(t),x}$ instead of $\tau(t,x(t),x)$. For technical reasons, we study separately the cases in which function $\tau_{t,x(t),\cdot}$ is nonincreasing and nondecreasing.

\subsection{Function $\tau(t,x(t),\cdot)$ is nonincreasing}

In this part we study problem (\ref{p1}) in the case that deviating function $\tau$ is nonincreasing in its third variable. We follow the spirit of \cite{jan} and we make the following definitions.

\begin{definition}\label{quasi} We say that $v,w \in {\cal C}(I_{\pm})$ are quasisolutions of (\ref{p1}) if $v_{|I_+}, w_{|I_+} \in AC(I_+)$ and
\begin{equation}\nonumber
\left\{
\begin{array}{ll}
v'(t)=f(t,v(t),w(\tau_{t,v(t),v}),v) \, \, \mbox{for a.a. $t \in I_+$,} \quad v(t)=\Lambda(v)+k(t) \, \, \mbox{for all $t \in I_-,$} \\
\\
w'(t)=f(t,w(t),v(\tau_{t,w(t),w}),w) \, \, \mbox{for a.a. $t \in I_+$,} \,\, w(t)=\Lambda(w)+k(t)\, \, \mbox{for all $t \in I_-.$}
\end{array}
\right.
\end{equation}
We say that $v_*,w^* \in Y \subset {\cal C}(I_{\pm})$ are the extremal quasisolutions of (\ref{p1}) in $Y$ if they are quasisolutions and $v_* \leq v$ and $w \leq w^*$ for any other couple of quasisolutions $v,w \in Y$ (in particular, $v_* \le w^*$ because $w^*$, $v_*$ are quasisolutions too).
\end{definition}

\begin{definition}\label{subsol} We say that $\alpha, \beta \in {\cal C}(I_{\pm})$ are, respectively, a lower and an upper solution for problem (\ref{p1}) if $\alpha_{|I_+}, \beta_{|I_+} \in AC(I_+)$ and the following inequalities hold:
\begin{equation}\nonumber
\left\{
\begin{array}{ll}
\alpha'(t) \leq f(t,\alpha(t),\beta(\tau_{t,\alpha(t),\alpha}),\alpha) \, \, \mbox{for a.a. $t \in I_+$,} \quad \alpha(t) \leq\Lambda(\alpha)+ k(t) \, \, \mbox{for all $t \in I_-,$} \\
\\
\beta'(t) \geq f(t,\beta(t),\alpha(\tau_{t,\beta(t),\beta}),\beta)\, \, \mbox{for a.a. $t \in I_+$,} \quad \beta(t) \geq \Lambda(\beta)+k(t)\, \, \mbox{for all $t \in I_-.$} \end{array}
\right.
\end{equation}

\end{definition}

Next we state and prove our main result on existence of extremal quasisolutions for problem (\ref{p1}) in case $\tau$ is nonincreasing in its third argument.

\begin{theorem}\label{main1} (Extremal quasisolutions) Suppose that either $r=0$ and $k$ is a fixed real number,
\noindent or $r>0$ and $k:I_- \longrightarrow \mathbb R$ is continuous on $I_-$ and nondecreasing on $[t_0-\hat r,t_0]$ for some $\hat r \in [0,r]$.

Assume
\medbreak
\begin{enumerate}
\item[($H_1$)]{(Lower and upper solutions) There exist $\alpha,\beta$,  lower and upper solution to (\ref{p1}), res\-pectively, such that $\alpha \leq \beta$ on $I_{\pm}$ and $\alpha$ and $\beta$ are nondecreasing on $[t_0-\hat r,t_0+L]$;}

\medbreak

   \item[$(H_2)$] ($L^1$ bounds) There exist $\psi_m, \psi_M \in L^1(I_+)$ such that for a.a. $t \in I_+$, all $x \in [\alpha(t),\beta(t)]$, and all $\gamma_i \in {\cal C}(I_{\pm})$ such that $\alpha \le \gamma_i \leq \beta$ on $I_{\pm}$ ($i=1,2$), we have
    $$0 \le \psi_m(t) \le f(t,x,\gamma_2(\tau(t,x,\gamma_1)),\gamma_1)\le \psi_M(t),$$
    \end{enumerate}

    \medbreak

\noindent
and let $
[\alpha,\beta]^+$ be the set of functions $\gamma \in {\cal C}(I_{\pm})$ such that $\alpha \leq \gamma \leq \beta$ on $I_{\pm}$, $\gamma_{|I_+} \in AC(I_+)$, $\psi_m \leq \gamma'\leq \psi_M$ a.e. on $I_+$, and $\gamma$ nondecreasing on $[t_0-\hat r,t_0+L]$.

\medbreak

 Finally, suppose that the following conditions hold:

 \medbreak
 \begin{enumerate}

\item[($H_3$)] (Discontinuous and non--monotone dependences with respect to $t$ and $x(t)$)
\begin{enumerate}
\item[(a)] (Measurability w.r.t. $t$) For all $x \in [\alpha(t_0),\beta(t_0+L)]$ and all $\gamma_i \in [\alpha,\beta]^+$ ($i=1,2$), the compositions

    $t \in \{ s \in I_+ \, : \, \alpha(s) \le x \le \beta(s)\} \longmapsto f(t,x,\gamma_2(\tau(t,x,\gamma_1)),\gamma_1)$, \\

    \noindent
    $t \in I_+ \longmapsto f(t,\alpha(t),\gamma_2(\tau(t,\alpha(t),\gamma_1)),\gamma_1)$, and \\

    \noindent $t \in I_+ \longmapsto f(t,\beta(t),\gamma_2(\tau(t,\beta(t),\gamma_1)),\gamma_1)$ are measurable;

    \medbreak

    \item[(b)] (Admissible discontinuities w.r.t. $x(t)$) For a.a. $t \in I_+$ and all $\gamma_i \in [\alpha,\beta]^+$ ($i=1,2$), the function $f(t,\cdot,\gamma_2(\tau(t,\cdot,\gamma_1)),\gamma_1)$ is continuous on $[\alpha(t),\beta(t)]
\backslash K(t)$, where $
K(t)= \cup_{n=1}^{\infty}K_n(t)$ and may depend on the choice of $\gamma_i$ ($i=1,2$),
and for
each $n \in \mathbb N$ and each $x \in K_n(t)$ we have
$$
\bigcap_{\varepsilon >0}
\overline{{\rm co}} f(t, x + \varepsilon B,\gamma_2(\tau_{t,x+\varepsilon B,\gamma_1}),\gamma_1)
 \bigcap DK_n(t,x)(1)   \subset
\{f(t,x,\gamma_2(\tau_{t,x,\gamma_1}),\gamma_1)\},
$$
where $B$ and $\overline{{\rm co}}$ are as in Lemma \ref{lemdis}.

    \end{enumerate}

    \medbreak

\item[($H_4$)] (Monotone functional dependences)

\begin{enumerate}
\item[(a)] For a.a. $t \in I_+$, for all $x \in [\alpha(t),\beta(t)]$, and for all $\gamma \in [\alpha,\beta]^+$ the function $f(t,x,\cdot,\gamma)$ is nonincreasing on $[\alpha(t_0-\hat r),\beta(t_0+L)]$; and for a.a. $t \in I_+$, for all $x \in [\alpha(t),\beta(t)]$, and for all $y \in [\alpha(t_0-\hat r),\beta(t_0+L)]$ the operator $f(t,x,y,\cdot)$ is nondecreasing on $[\alpha,\beta]^+$, i.e., for $\gamma_i \in [\alpha,\beta]^+$ ($i=1,2$) the relation $\gamma_1 \le \gamma_2$ on $I_{\pm}$ implies $f(t,x,y,\gamma_1) \le f(t,x,y,\gamma_2)$;

\medbreak

\item[(b)] For a.a. $t \in I_+$ and all $x \in [\alpha(t),\beta(t)]$, the function

\noindent $\tau(t,x,\cdot):[\alpha,\beta]^+ \longrightarrow [t_0-\hat r,t_0+L]$ is nonincreasing, i.e., for $\gamma_i \in [\alpha,\beta]^+$ ($i=1,2$) the relation $\gamma_1 \le \gamma_2$ on $I_{\pm}$ implies $$t_0-\hat r \le \tau(t,x,\gamma_2) \le \tau(t,x,\gamma_1) \le t_0+L;$$

\medbreak

\item[(c)] $\Lambda$ is nondecreasing on $[\alpha,\beta]^+$.
\end{enumerate}

\end{enumerate}
\medbreak
Then problem (\ref{p1}) has the extremal quasisolutions $v_*,w^*$ in $[\alpha,\beta]^+,$  which satisfy \begin{equation}\label{carac}
(v_*,w^*)=\min_{\preceq} \{(v,w) \, : \,  \mbox{$w,v$ are coupled lower and upper solutions in $[\alpha,\beta]^+$ for } (\ref{p1})\},
\end{equation}
where $\preceq$ is the partial ordering defined in the proof of Theorem \ref{multi} for $X=\mathcal{C}(I_{\pm})$.
\end{theorem}

\noindent {\bf Proof.}  We reduce our problem to that of finding the extremal coupled fixed points of an adequate multivalued operator. To do so, and in order to apply Theorem \ref{multi}, we consider the space $X=\mathcal{C}(I_{\pm})$ endowed with its usual metric and pointwise ordering. In
the subset
\begin{align*}
Y=\left\{ \gamma \in X \, : \,  \right.  & \mbox{$\gamma_{|I_+} \in AC(I_+)$, $\psi_m \leq \gamma'\leq \psi_M$ a.e. on $I_+$,} \\
& \left. \mbox{and $\gamma$ nondecreasing on $[t_0-\hat r,t_0+L]$} \right\}
\end{align*}
 we have the ordered interval $[\alpha,\beta]^+$ introduced in the statement.

We define a multivalued operator ${\cal A}:[\alpha,\beta]^+ \times [\alpha,\beta]^+ \longrightarrow 2^{[\alpha,\beta]^+} \setminus \emptyset$ as follows: for each pair $(\gamma_1,\gamma_2) \in [\alpha,\beta]^+ \times [\alpha,\beta]^+$ we define ${\cal A}(\gamma_1,\gamma_2)$ as the set of solutions in $[\alpha,\beta]^+$ to the initial value problem
\begin{equation}\label{ivp}
z'(t)=f(t,z(t),\gamma_2(\tau_{t,z(t),\gamma_1}),\gamma_1) \mbox{ a.e. on $I_+$}, \quad z=\Lambda({\gamma_1})+k \, \, \mbox{on $I_-.$}
\end{equation}
Lemma \ref{lemdis} and assumptions $(H_1)$ through $(H_3)$ guarantee that ${\cal A}(\gamma_1,\gamma_2) \neq \emptyset$ for all $(\gamma_1,\gamma_2) \in [\alpha,\beta]^+ \times [\alpha,\beta]^+$ and that, moreover, the corresponding single--valued operators $A_{\pm}$ (in the terms of Theorem \ref{multi}) are well--defined.

\medbreak
\noindent
{\it Claim 1 -- The operators $A_{\pm}$ are mixed monotone.} Let $\gamma_1,\overline{\gamma_1},\gamma_2 \in [\alpha,\beta]^+$ be such that $\gamma_1 \leq \overline{\gamma_1}$, and put $\xi = A_+(\gamma_1,\gamma_2)$, $\overline{\xi}=A_+(\overline{\gamma_1},\gamma_2)$. Condition $(H_4)-(c)$ ensures that $\xi \leq \Lambda(\overline{\gamma_1})+k$ on $I_-$, and, on the other hand, $(H_4)-(a)$, the monotonicity of $\gamma_2$ on $[t_0-\hat{r},t_0+L]$, and $(H_4)-(b)$ guarantee that
\begin{equation}\label{des1}
\xi'(t)=f(t,\xi(t),\gamma_2(\tau_{t,\xi(t),\gamma_1}),\gamma_1) \leq f(t,\xi(t),\gamma_2(\tau_{t,\xi(t),\overline{\gamma_1}}),\overline{\gamma_1}) \quad \mbox{a.e. on $I_+$},
\end{equation}
so $\xi$ is a lower solution for problem (\ref{ivp}) with $\gamma_1$ replaced by $\overline{\gamma_1}$, whose greatest solution in $[\alpha,\beta]^+$ is $\overline{\xi}$. Hence $\xi \leq \overline{\xi}$ on $I_{\pm}$, by virtue of Lemma \ref{lemdis}, thus proving that $A_+(\cdot,\gamma_2)$ is nondecreasing. Similar arguments prove that $A_-(\cdot, \gamma_2)$ is nondecreasing and that $A_{\pm}(\gamma_1,\cdot)$ are nonincreasing.

\medbreak
\noindent
{\it Claim 2 -- The operators $A_{\pm}$ satisfy the m.m.c.p.} Let $(v_j)_{j=1}^{\infty}$ and $(w_j)_{j=1}^{\infty}$ be sequences in $[\alpha,\beta]^+$, one being nondecreasing and the other nonincreasing. The mixed monotoni\-ci\-ty property implies that $(z_j)_{j=1}^{\infty}=(A_{+(-)}(v_j,w_j))_{j=1}^{\infty}$ is a monotone sequence, which is bounded too, so its pointwise limit, which we denote by $z$, exists. It remains to prove that $z \in Y$ and $(z_j)_{j=1}^{\infty}$ converges uniformly on $I_{\pm}$. For $s,t \in I_+$, $s<t$, and $j \in \mathbb N$, condition $(H_2)$ yields
$$
\int_s^t \psi_m(r) \, dr \le  z_j(t)-z_j(s) = \int_s^t f(r,z_j(r),w_j(\tau_{r,z_j(r),v_j}),v_j) \, dr \leq \int_s^t \psi_M(r) \, dr,
$$
and going to the limit as $j$ tends to infinity we obtain
$$
\int_s^t \psi_m(r) \, dr\le z(t)-z(s) \leq \int_s^t{\psi_M(r)\, dr},$$
which implies that $z_{|I_+} \in AC(I_+)$ and $\psi_m \leq z'\leq \psi_M$ a.e. on $I_+$. \\
On the other hand, for all $t \in I_-$ we have $$z(t)=\lim_{j \to \infty}z_j(t)=\lim_{j \to \infty}\Lambda({v_j})+k(t),$$
   so $z_{|I_-}$ is continuous on $I_-$ and nondecreasing on $[t_0-\hat r,t_0]$. Therefore $z \in Y$. Finally, the monotone convergence of $(z_j)_{j=1}^{\infty}$ to the continuous function $z$ implies that this convergence is uniform on $I_{\pm}$, by virtue of Dini's Theorem. Claim 2 is proven.

   \medbreak

Theorem \ref{multi} guarantees now that ${\cal A}$ has the extremal coupled fixed points $v_*, w^*$, which are quasisolutions of (\ref{p1}) in $[\alpha,\beta]^+$. Moreover, if $v$ and $w$ are quasisolutions of (\ref{p1}) in $[\alpha,\beta]^+$, then $(v,w)$ satisfies $(CFP)$, so we have $v_*\le v$ and $w\le w^*$ by virtue of $(ECFP)$. Therefore $v_*$ and $w^*$ are the extremal quasisolutions of (\ref{p1}) in $[\alpha,\beta]^+$.

\medbreak

Finally, characterization (\ref{minmax}) ensures that
\begin{equation}\label{carac2}
(v_*,w^*)=\min_{\preceq} \{(v,w) \in [\alpha,\beta]^+ \times [\alpha,\beta]^+ \, : \, A_-(v,w) \leq v, \ w \leq A_+(w,v)\},
\end{equation}
so our last purpose is to check that both expressions (\ref{carac}) and (\ref{carac2}) are equivalent. \\
Put
$$
\begin{array}{ll}
S&=\{(v,w) \, : \, w,v \mbox{ are coupled lower and upper solutions in $[\alpha, \beta]^+$ for } (\ref{p1})\}, \\
&
\\
\overline{S}&=\{(v,w) \in [\alpha,\beta]^+ \times [\alpha,\beta]^+ \, : \, A_-(v,w) \leq v, \ w \leq A_+(w,v)\}.
\end{array}$$

Let $(v,w) \in S$ and denote $x=A_-(v,w)$, which is the least solution in $[\alpha,\beta]^+$ for the problem
\begin{equation}\label{procarac1}
x'(t)=f(t,x(t),w(\tau_{t,x(t),v}),v) \mbox{ on $I_+,$ } \quad x=\Lambda(v) + k \mbox{ on $I_-$},
\end{equation}
so (\ref{least}) implies that $x \leq v$ because $v$ is an upper solution for (\ref{procarac1}). In an analogous way, we check that $w \leq A_+(w,v)$ and therefore $S \subset \overline{S}$, so $\min_{\preceq} \overline{S} \preceq \min_{\preceq} S.$ \\
On the other hand, given $(\overline{v},\overline{w}) \in \overline{S}$ functions $v=A_-(\overline{v},\overline{w})$ and ${w}=A_+(\overline{w},\overline{v})$ satisfy that $(v,v) \preceq (\overline{v},\overline{w})$ and by condition $(H_4)$ they are, respectively, an upper and a lower solution in $[\alpha,\beta]^+$ for problem (\ref{p1}). Therefore, given an element $(\overline{v},\overline{w}) \in \overline{S}$ there exists $(v,w) \in S$ such that $(v,w) \preceq (\overline{v},\overline{w})$ and then $\min_{\preceq} S \preceq \min_{\preceq} \overline{S},$ and this ends the proof.
\qed

Now we have sufficient conditions for the existence of the extremal quasisolutions to (\ref{p1}) we will show that, under some more assumptions, both functions are identical, thus defining a solution. To achive this goal, we employ the following maximum principle, which generalizes \cite[Lemma 1]{jan}:

\begin{lemma}
\label{lem1}
Let $p \in {\cal C}(I_{\pm})$, $\Psi \in L^1(I_+,[0,\infty))$, and $\lambda \in [0,1)$.

If $p_{|I_+} \in AC(I_+)$ and
\begin{eqnarray}
\label{el1}
p'(t) \le \Psi(t) \max_{s \in I_{\pm}}p(s) \quad \mbox{for a.a. $t \in I_+$,} \\
\label{el2}
p(t)\le \lambda \max_{s \in I_{\pm}}p(s)\quad \mbox{for all $t \in I_-$,}\\
\nonumber
\mbox{and}\\
\label{el3}
\lambda+\int_{t_0}^{t_0+L}{\Psi(s)ds}<1,
\end{eqnarray}
then $p(t) \le 0$ for all $t \in I_{\pm}$.
\end{lemma}

\noindent
\noindent {\bf Proof.}  Let $t_1\in I_{\pm}$ be such that
$$p(t_1)=\max_{s \in I_{\pm}}p(s),$$
and assume, reasoning by contradiction, that $p(t_1)>0.$

On the one hand, if $t_1 \le t_0$ then condition (\ref{el2}) yields
$$p(t_1) \le \lambda p(t_1)<p(t_1),$$
a contradiction.

On the other hand, if $t_1>t_0$ then we integrate in (\ref{el1}) between $t_0$ and $t_1$, we use (\ref{el2}) for $t=t_0$, and we obtain
$$p(t_1) \le p(t_0)+p(t_1)\int_{t_0}^{t_0+L}{\Psi(s)ds} \le p(t_1)\left( \lambda+\int_{t_0}^{t_0+L}{\Psi(s)ds}\right),$$
once again a contradiction, by virtue of (\ref{el3}).
\qed

Next we give sufficient conditions for (\ref{p1}) to have a unique solution between given lower and upper solutions.

\begin{theorem} (Unique solutions) \label{main2} Suppose that all the assumptions in Theorem \ref{main1} hold, and assume that there exists $\hat \psi \in L^1(t_0-\hat r,t_0)$ such that for $s,t \in [t_0-\hat r,t_0]$, $s \le t$, we have
\begin{equation}
\label{bv-}
k(t)-k(s) \le \int_s^t{\hat \psi(r) dr}.
\end{equation}
Moreover, we assume that the following group of conditions is satisfied:

\medbreak

\begin{enumerate}

    \item[($H_5$)] (One--sided Lipschitz functional dependences)
    \begin{enumerate}
    \item[(a)] Let $L_1, L_2 \in L^1(I_+,[0,\infty))$ be such that for a.a. $t \in I_+$, for all $x \in [\alpha(t),\beta(t)]$, and for all $\gamma_i \in [\alpha,\beta]^+$ ($i=1,2$), the relations $\alpha(\tau_{t,x,\beta})\le u \le v \le \beta(\tau_{t,x,\alpha})$ and $\gamma_1 \le \gamma_2$ on $I_{\pm}$ imply
  \begin{equation}\label{h4}
f(t,x,u,\gamma_2)-f(t,x,v,\gamma_1) \le L_1(t)(v-u)+L_2(t)\max_{s \in I_{\pm}}(\gamma_2(s)-\gamma_1(s));
\end{equation}

\medbreak

\item[(b)] Let $L_3 \in L^1(I_+,[0,\infty))$ be such that $L_1 \cdot L_3 \in L^1(I_+,[0,\infty))$ and for a.a. $t \in I_+$, for all $x \in [\alpha(t),\beta(t)]$, for $\gamma_i \in [\alpha,\beta]^+$ ($i=1,2$), $\gamma_1 \le \gamma_2$, we have
    \begin{equation}
    \label{lipbet}
    \int_{\tau_{t,x,\gamma_2}}^{\tau_{t,x,\gamma_1}}{\tilde \psi (s)ds} \le L_3(t)
    \max_{s \in I_{\pm}}(\gamma_2(s)-\gamma_1(s)),
    \end{equation}
    where $\tilde \psi=\hat \psi$ a.e. on $[t_0-\hat r,t_0]$ ($\hat{\psi}$ as in (\ref{bv-}))$, \tilde \psi=\psi_M$ a.e. on $I_+$ ($\psi_M$ as in $(H_2)$);

    \medbreak

    \item[(c)] There exists $\lambda \in [0,1)$ such that for $\gamma_i \in [\alpha,\beta]^+$ ($i=1,2$) with $\gamma_1 \le \gamma_2$, we have
    $$\Lambda({\gamma_2})-\Lambda({\gamma_1}) \le \lambda
    \max_{s \in I_{\pm}}(\gamma_2(s)-\gamma_1(s)).$$
    \end{enumerate}

    \medbreak

    \item[$(H_6)$] (One--sided Lipschitz condition w.r.t. $x(t)$)

    There exists $L_4 \in L^1(I_+,[0,\infty))$ such that for a.a. $t \in I_+$ and all $\gamma_i \in [\alpha,\beta]^+$ with $\gamma_i=c_i+k$ on $I_-$ for some constants $c_i \in \mathbb R$ ($i=1,2$) the relation $\alpha(t) \le x \le y \le \beta(t)$ implies
  \begin{equation}\label{h41}
f(t,y,\gamma_1(\tau(t,y,\gamma_2)),\gamma_2)-f(t,x,\gamma_1(\tau(t,x,\gamma_2)),\gamma_2) \le L_4(t)(y-x).
\end{equation}
\end{enumerate}
Then problem (\ref{p1}) has a unique solution in $[\alpha, \beta]^+$ provided that  \begin{equation}
\label{small}
\lambda + \int_{t_0}^{t_0+L}{L_1(s)ds}+\int_{t_0}^{t_0+L}{L_2(s)ds}+ \int_{t_0}^{t_0+L}{L_1(s)L_3(s)ds} +\int_{t_0}^{t_0+L}{L_4(s)ds}<1.
\end{equation}
\end{theorem}
\noindent
\noindent {\bf Proof.}  By Theorem \ref{main1} we know that problem (\ref{p1}) has the extremal quasisolutions $v_*$, $w^*$ in $[\alpha, \beta]^+$ and $v_* \leq w^*$ on $I_{\pm}$. Let us prove that $w^* \le v_*$ on $I_{\pm}$; to do so, we
 define $p=w^*-v_*$ and we use our assumptions over $f$ and $\tau$ to establish the following inequalities for a.a. $t \in I_+$:
\begin{align*}
p'(t)&=f(t,w^*(t),v_*(\tau_{t,w^*(t),w^*}),w^*)-f(t,v_*(t),v_*(\tau_{t,v_*(t),w^*}),w^*) \\
& \mbox{\quad}+f(t,v_*(t),v_*(\tau_{t,v_*(t),w^*}),w^*)-f(t,v_*(t),w^*(\tau_{t,v_*(t),v_*}),v_*) \\
&\le L_4(t)(w^*(t)-v_*(t))+L_1(t)(w^*(\tau_{t,v_*(t),v_*})-v_*(\tau_{t,v_*(t),w^*}))+L_2(t)\max_{s \in I_{\pm}}p(s)\\
& \le (L_2(t)+L_4(t))\max_{s \in I_{\pm}}p(s)+L_1(t)p(\tau_{t,v_*(t),v_*})+L_1(t)(v_*(\tau_{t,v_*(t),v_*})-v_*(\tau_{t,v_*(t),w^*})) \\
& \le (L_1(t)+L_2(t)+L_4(t))\max_{s \in I_{\pm}}p(s)+L_1(t)\int_{\tau_{t,v_*(t),w^*}}^{\tau_{t,v_*(t),v_*}}{\tilde \psi(s)ds}\\
 & \le (L_1(t)+L_2(t)+L_1(t)L_3(t) + L_4(t))\max_{s \in I_{\pm}}p(s).
\end{align*}
On the other hand, for all $t \in I_-$ we have
$$p(t)=\Lambda({w^*})-\Lambda({v_*}) \le \lambda \, \max_{s \in I_{\pm}}p(s),$$
thus Lemma \ref{lem1} and condition (\ref{small}) imply $p \le 0$ on $I_{\pm}$. Hence, $v_*=w^*$ on $I_{\pm}$, and therefore $v_*$ is a solution of (\ref{p1}).

Finally, if $y \in [\alpha,\beta]^+$ is another solution of (\ref{p1}) then $y, \, y$ are quasisolutions, and $(ECFP)$ with $P=[\alpha,\beta]^+$ implies $v_* \leq y \leq w^*=v_*$, so $y=v_*$.
\qed

\begin{example}\label{ej1} First, let $\{q_m\}_{m=1}^{\infty}$ be an enumeration of all rational numbers in $[0,\infty)$ such that $q_1=4$ and consider the function
 $$\phi_1(x)=1-\sum_{\{m \,:\, q_m < x\}} 2^{-m} \quad (x \in [0,\infty)),$$
  which is decreasing, continuous exactly on the set of irrational numbers in $[0,\infty)$, and sa\-tis\-fies $0 < \phi_1 \leq 1$ on $[0, \infty)$. Notice that the choice $q_1=4$ implies that $\phi_1(5) \geq 1/2$, which will be used in our argumentations.

\medbreak

  Second, let $I_+=[0,1]$ and fix $\rho \in (0,1)$. For each $m \in \mathbb N$ and $n \in \mathbb Z \backslash \{0\}$ we define the sets $$K_{m,n}(t)=\left\{3 \dfrac{n}{|n|}t + \rho m - \dfrac{1}{|n|}\right\} \quad (t \in I_+),$$
 and for each $(t,x) \in I_+ \times [0,\infty)$ let
 $$\phi_2(t,x)=
\left\{
\begin{array}{cl}
1, \ &\mbox{if } x \in K(t)=\bigcup_{m \in \mathbb N} \, \bigcup_{n \in \mathbb Z \backslash \{0\}} \, K_{m,n}(t), \\
\\
 \left( \dfrac{1}{5}tx \right)^{\mu} \, \cos x, \ &\mbox{if } x \notin K(t),
\end{array}
\right.$$
where $\mu$ is a positive constant. Notice that for any $t \in I_+$ the function $\phi_2(t,\cdot)$ is continuous on $[0,\infty) \setminus K(t)$, but it is neither monotone nor bounded on $[0,\infty)$ if $t>0$. Notice also that $x =\pm 3\, t+ \rho m$ ($m \in \mathbb N$) are accumulation points of $K(t)$.

\bigbreak

We are concerned with the functional problem
\begin{equation}\label{ex1}
\left\{
\begin{array}{lcl}
x'(t)&=&\sigma_1(t) \, \phi_1(x(t)) + \sigma_2(t) \, \phi_2(t,x(t)) - \varepsilon(t)   \, x \left( 
\frac{1-x(0)\, |\sin x(t)|\, t}{2} \right), \ t \in I_+, \\
\\
x(t)&=&\int_{I_+}\omega(s)x(s) \, ds + \cos t, \ t \in I_-=[-1,0].
\end{array}
\right.
\end{equation}
Problem (\ref{ex1}) is the particular case of (\ref{p1}) corresponding to
$$\tau_{t,x,\gamma}=\frac{1-\gamma(0)\, |\sin x |\, t}{2} \quad \mbox{for all $(t,x,\gamma)
\in I_+ \times \mathbb R \times {\cal C}(I_{\pm})$,}$$
$$f(t,x,y)=\sigma_1(t) \, \phi_1(x) + \sigma_2(t) \, \phi_2(t,x) - \varepsilon(t)   \, y \quad \mbox{for all $(t,x,y)
\in I_+ \times \mathbb R^2$,}$$
$\Lambda(\gamma)=\int_{I_+}\omega(s)\gamma(s)$ for all $\gamma \in {\cal C}(I_{\pm})$, and $k(t)= \cos t$ for all $t \in I_-$.

\medbreak

 We assume that $\varepsilon, \, \sigma_i \in L^{\infty}(I_+,[0,\infty))$ ($i=1,2$), $\omega$ is an integrable nonnegative (a.e.) function, there exist $\delta_1 >0$ and $\delta_2\in (0,3)$ such that $\|\sigma_1+\sigma_2\|_{\infty} \leq 3-\delta_2$, for a.a. $t \in I_+$ we have \begin{equation}\label{epsilon}  \sigma_1(t) \geq \left( 2 + \frac{3-\delta_2}{2} \right) \varepsilon(t),\end{equation} \begin{equation}\label{sigma}  \sigma_1(t) \geq 2 \left( \delta_1 + \sigma_2(t) \left( \dfrac{5-\delta_2}{5}\right)^{\mu} + \varepsilon(t) (5-\delta_2) \right), \end{equation} and \begin{equation}\label{omega} \|\omega\|_1 \leq (5-\delta_2)^{-1}. \end{equation}

We consider now functions $\alpha(t)=0$ and $\beta(t)=2 + \chi_{[0,1]} (t)(3-\delta_2)t$, $t \in I_{\pm},$ and we will check that, under certain assumptions, $\alpha$ and $\beta$ are, respectively, a lower and an upper solution for problem (\ref{ex1}). Indeed, for $t \in I_+$ we have 
$$
f(t,\alpha(t),\beta(\tau_{t,\alpha(t),\alpha}))=\sigma_1(t)-\varepsilon(t)\left( 2 + \frac{3-\delta_2}{2} \right),$$
so by (\ref{epsilon}) we obtain that $\alpha'(t) \leq f(t,\alpha(t),\beta(\tau_{t,\alpha(t),\alpha}))$ for $t \in I_+$. Moreover, for $t \in I_-$ we have $0=\alpha(t) \leq \cos t$. \\
On the other hand, for a.a. $t \in I_+$ it is $\phi_1(\beta(t)) \leq 1$, $\phi_2(t,\beta(t)) \leq 1$, and then
$$
f(t,\beta(t),\alpha(\tau_{t,\beta(t),\beta})) \leq \sigma_1(t)+\sigma_2(t) \leq 3-\delta_2 = \beta'(t).$$
Finally, for $t \in I_-$ we have
$$
\int_{I_+} \omega(s) \beta(s) \, ds + \cos t \leq (5-\delta_2) \|\omega\|_1 + \cos t \leq 2=\beta(t),$$
where the previous inequality is fulfilled by virtue of (\ref{omega}). 

\bigbreak

Now we will show that this problem has extremal quasisolutions between $\alpha$ and $\beta$.

Notice that for a.a. $t \in I_+$ and $x \leq \beta(t)$ we have $x \le 5$, so $\phi_1(x) \geq \phi_1(5) \geq 1/2$, and
$$\phi_2(t,x) \geq - \left( \dfrac{5-\delta_2}{5}\right)^{\mu},$$
hence for a.a. $t \in I_+$, $0 \leq x \leq \beta(t)$ and $0\leq y \leq 5-\delta_2$ we have $f(t,x,y) \geq \delta_1$ by virtue of (\ref{sigma}). Therefore, condition $(H_2)$ in Theorem \ref{main1} is satisfied (take for example $\psi_m = \delta_1$ and $\psi_M=\sigma_1 + \sigma_2$). The  monotonicity conditions $(H_4)$ are also satisfied, so we only have to check conditions $(H_3)-(a)$ and $(H_3)-(b)$. \\

In order to check measurability conditions we use proposition 3.2. in \cite{cp}, a result on {\it superpositional measurability}, that is, a result that provides us with sufficient conditions to guarantee that the composition operator associated to a function $g:I \times \mathbb R \rightarrow \mathbb R$ maps continuous functions into measurable ones.  

So, let $\xi_i$ ($i=1,2,3$) be absolutely continuous functions on $I_+$, with $\xi_1 \geq 0$. The mapping $$t \in \{s \in I_+ \, : \, \alpha(s) \le \xi_1(s) \le \beta(s) \} \longmapsto f(t,\xi_1(t),\xi_2(\tau(t,\xi_1(t),\xi_3)))$$ can be written as
$$
f(\cdot,\xi_1(\cdot),\xi_2(\tau(\cdot,\xi_1(\cdot),\xi_3)))=\sigma_1(\cdot) \widetilde{\phi_1}(\cdot) + \sigma_2(\cdot) \widetilde{\phi_2}(\cdot) - \varepsilon(\cdot) \widetilde{\phi_3}(\cdot),$$
where $\sigma_1, \sigma_2$ and $\varepsilon$ are measurable by assumption and:
\begin{enumerate}
\item $\widetilde{\phi_1}=\phi_1 \circ \xi_1$ is monotone, so measurable,
\item $\phi_2$ satisfies the conditions of proposition 3.2 in \cite{cp}, so $\widetilde{\phi_2}(\cdot)=\phi_2(\cdot,\xi_1(\cdot))$ is measurable,
\item $\widetilde{\phi_3}(\cdot)=\xi_2(\tau(\cdot,\xi_1(\cdot),\xi_3))$ is continuous, so measurable.
\end{enumerate}
Then we obtain that required measurability conditions $(H_3)-(a)$ are satisfied.\\

Finally, discontinuities with respect to $x(t)$ are due to functions $\phi_1$ and $\phi_2$. In the first case the discontinuity set can be written as a countable union of horizontal lines and in the second case the discontinuity set can also be written as a countable union of lines, now with derivative $\pm 3$. So, condition $(H_3)-(b)$ in Theorem \ref{main1} holds because between $\alpha$ and $\beta$ we have $0 < \delta_1 \leq f \leq 3- \delta_2 <3$. \\

In conclusion, by application of Theorem \ref{main1} we obtain that problem (\ref{ex1}) has the extremal quasisolutions between $\alpha$ and $\beta$. \\

Unfortunately, we cannot use Theorem \ref{main2} in order to look for an unique solution between $\alpha$ and $\beta$ for problem (\ref{ex1}) because function $\phi_2$ doesn't sa\-tis\-fy the monotone conditions that hypothesis $(H_6)$ requires. However, we will show that if we drop this term (that is, if we put $\sigma_2 \equiv 0$) then we can apply Theorem \ref{main2}. \\

First of all, for a.a. $t \in I_+$ and all $x \in [\alpha(t),\beta(t)]$ we have
$$
f(t,x,u)-f(t,x,v)=\varepsilon(t)(v-u)$$
whenever $\alpha(\tau_{t,x,\beta}) \leq u \leq v \leq \beta(\tau_{t,x,\alpha})$, so condition $(H_5)-(a)$ holds with $L_1=\varepsilon$ and $L_2=0$.   

Second, for a.a. $t \in I_+$, all $x \in [\alpha(t),\beta(t)],$ and all $\gamma_1,\gamma_2 \in [\alpha,\beta]^+$ with $\gamma_1 \leq \gamma_2$ we have  
$$
\int_{\tau_{t,x,\gamma_2}}^{\tau_{t,x,\gamma_1}} \widetilde{\psi}(s) \, ds \leq \dfrac{t \, \max\{1,\|\sigma_1\|_{\infty}\}}{2} \, (\gamma_2(0)-\gamma_1(0)),$$
where $\widetilde{\psi}(t)=-\sin t, \, t \in I_-, $ and $\widetilde{\psi}=\sigma_1 $ on $I_+$.
Therefore condition $(H_5)-(b)$ holds with $L_3(t)=t \, \max\{1,\|\sigma_1\|_{\infty}\}/2$.

Third, for $\gamma_1,\gamma_2 \in [\alpha,\beta]^+$ with $\gamma_1 \leq \gamma_2$ we have
$$
\Lambda(\gamma_{2}) - \Lambda(\gamma_{1}) \leq \|\omega\|_1 \, \max_{s \in I_{\pm}} (\gamma_2(s)-\gamma_1(s)),$$
so condition $(H_5)-(c)$ hold with $\lambda=\|\omega\|_1$. \\

Finally, for a.a. $t \in I_+$, for $\gamma_i \in [\alpha,\beta]^+$ ($i=1,2$), and $\alpha(t) \le x \le y \le \beta(t)$ we have
\begin{align*}
f(t,y,\gamma_1(\tau_{t,y,\gamma_2}))-f(t,x,\gamma_1(\tau_{t,x,\gamma_2})) & \le \varepsilon(t)(\gamma_1(\tau_{t,x,\gamma_2})-\gamma_1(\tau_{t,y,\gamma_2})) \\
& \le \varepsilon(t)\| \sigma_1Ê\|_{\infty}\dfrac{1}{2} \gamma_2(0)t \left| |\sin x|-|\sin y|\right| \\
& \le \varepsilon(t)\| \sigma_1Ê\|_{\infty} \, t (y-x),
\end{align*}
so $(H_6)$ is satisfied with $L_4(t)=\varepsilon(t)\| \sigma_1Ê\|_{\infty} \, t$, $t \in I_+$.

\medbreak

Consequently, problem (\ref{ex1}) has a unique solution in $[\alpha,\beta]^+$ whenever the previously defined $L_i$ ($i=1,2,3,4$) and $\lambda$ satisfy (\ref{small}), which holds, in particular, if
$$
\|\omega\|_1 + \|\varepsilon\|_{\infty}\left(1 + \dfrac{3\max\{1,\|\sigma_1 \|_{\infty}\}}{4} \right) < 1.$$

\end{example}

\subsection{Function $\tau(t,x(t),\cdot)$ is nondecreasing}
Changing the monotonicity of $\tau$ in its third variable requires somes changes in the definitions of quasisolutions and lower and upper solutions, but then we have existence results which are similar to those in the previous section.
\begin{definition} \label{quasi2}We say that $v,w \in {\cal C}(I_{\pm})$ are quasisolutions of (\ref{p1}) if $v_{|I_+}, w_{|I_+} \in AC(I_+)$ and
\begin{equation}\nonumber
\left\{
\begin{array}{ll}
v'(t)=f(t,v(t),w(\tau_{t,v(t),w}),v) \, \, \mbox{for a.a. $t \in I_+$,} \quad v(t)=\Lambda(v)+k(t) \, \, \mbox{for all $t \in I_-,$} \\
\\
w'(t)=f(t,w(t),v(\tau_{t,w(t),v}),w) \, \, \mbox{for a.a. $t \in I_+$,} \,\, w(t)=\Lambda(w)+k(t)\, \, \mbox{for all $t \in I_-.$}
\end{array}
\right.
\end{equation}
We say that $v_*,w^* \in Y \subset {\cal C}(I_{\pm})$ are the extremal quasisolutions of (\ref{p1}) in $Y$ if they are quasisolutions and $v_* \leq v$ and $w \leq w^*$ for any other couple of quasisolutions $v,w \in Y$.
\end{definition}

\begin{definition}
 \label{subsol2}
 We say that $\alpha, \beta \in {\cal C}(I_{\pm})$ are, respectively, a lower and an upper solution for problem (\ref{p1}) if $\alpha_{|I_+}, \beta_{|I_+} \in AC(I_+)$ and the following inequalities hold:
\begin{equation}\nonumber
\left\{
\begin{array}{ll}
\alpha'(t) \leq f(t,\alpha(t),\beta(\tau_{t,\alpha(t),\beta}),\alpha) \, \, \mbox{for a.a. $t \in I_+$,} \quad \alpha(t) \leq\Lambda(\alpha)+ k(t) \, \, \mbox{for all $t \in I_-,$} \\
\\
\beta'(t) \geq f(t,\beta(t),\alpha(\tau_{t,\beta(t),\alpha}),\beta)\, \, \mbox{for a.a. $t \in I_+$,} \quad \beta(t) \geq \Lambda(\beta)+k(t)\, \, \mbox{for all $t \in I_-.$} \end{array}
\right.
\end{equation}
\end{definition}
Our next result is the analogue to Theorem \ref{main1} in case $\tau$ is nondecreasing in its third argument.

\begin{theorem}\label{main3} Suppose that all the conditions in Theorem \ref{main1} are fulfilled with lower and upper solutions in the sense of Definition \ref{subsol2} and assumption $(H_4)-(b)$ replaced by
\begin{enumerate}

\item[($H_4$)]
\begin{enumerate}
\item[(b')] For a.a. $t \in I_+$ and all $x \in [\alpha(t),\beta(t)]$, the function

\noindent $\tau(t,x,\cdot):[\alpha,\beta]^+ \longrightarrow [t_0-\hat r,t_0+L]$ is nondecreasing, i.e., for $\gamma_i \in [\alpha,\beta]^+$ ($i=1,2$) the relation $\gamma_1(t) \le \gamma_2(t)$ for all $t \in I_{\pm}$ implies $$t_0-\hat r \le \tau(t,x,\gamma_1) \le \tau(t,x,\gamma_2) \le t_0+L.$$
\end{enumerate}
\end{enumerate}
\medbreak
Then problem (\ref{p1}) has the extremal quasisolutions (in the sense of Definition \ref{quasi2}) in $[\alpha,\beta]^+$, and they satisfy (\ref{carac}).
\end{theorem}

\noindent {\bf Proof.}  Repeat the proof of Theorem \ref{main1} with the following multivalued operator: for each pair $(\gamma_1,\gamma_2) \in [\alpha,\beta]^+ \times [\alpha,\beta]^+$ let ${\cal A}(\gamma_1,\gamma_2)$ be the set of solutions in $[\alpha,\beta]^+$ to the initial value problem
\begin{equation}\label{ivp2}
z'(t)=f(t,z(t),\gamma_2(\tau_{t,z(t),\gamma_2}),\gamma_1) \mbox{ a.e. on $I_+$}, \,\, z(t)=\Lambda({\gamma_1})+k(t) \, \, \mbox{on $I_-.$}
\end{equation}
\qed

As a consequence of Theorem \ref{main3} and Lemma \ref{lem1} we have the following uniqueness result.

\begin{theorem}\label{main4} In the conditions of Theorem \ref{main3}, suppose that conditions $(H_5)$ and $(H_6)$ in Theorem \ref{main2} hold with inequality (\ref{lipbet}) replaced by
 \begin{equation}
    \label{lipbet2}
    \int_{\tau_{t,x,\gamma_1}}^{\tau_{t,x,\gamma_2}}{\tilde \psi(s)ds} \le L_3(t)
    \max_{s \in I_{\pm}}(\gamma_2(s)-\gamma_1(s)).
    \end{equation}
If (\ref{small}) is fulfilled, then problem (\ref{p1}) has a unique solution in the functional interval $[\alpha,\beta]^+$.
\end{theorem}

\subsection{Discussion about the assumptions}
First, the constant $\hat r \in [0, r]$ in the assumptions specifies the range of $\tau$, where solutions must be monotone. We introduce $\hat r$ in order to relax our assumptions in the sense that they must be checked just in a part of $I_{\pm}$.
 
 \bigbreak
 
\noindent
$(H_1)$ Existence of lower and upper solutions is surely the assumption which is more difficult to check in practice. That is why it is useful to have sufficient conditions which imply $(H_1)$, and we include some in Section 4.

\bigbreak

\noindent
$(H_2)$ This type of integrable bounds is commonplace in this kind of research, and not very stringent: notice that $(H_2)$ is immediately fulfilled if nonlinearities are continuous (and nonnegative), thanks to boundedness on compact sets.

Notice also that the closer $\psi_m$ and $\psi_M$ can be chosen, then the smaller the set  $[\alpha,\beta]^+$ is, which might help checking the subsequent assumptions. This is the only reason for introducing a lower bound $\psi_m$.

Whether our previous results remain valid with locally integrable bounds, in the spirit of \cite{bp}, or, more important, with sign--changing right--hand sides, are open problems. We partially solve the second problem in Section 5, where we show that $f$ can change sign in case $\tau$ depends only on $(t,x(t))$.

\bigbreak

\noindent
$(H_3)-(a)$ Measurability of compositions is a very mild assumption, yet it is not always easy to check. Of course, we have it for granted with continuous or, more generally, Carath\'eodory nonlinearities, but it has to be explicitly required in our case because neither $f$ nor $\tau$ need be continuous with respect to $x(t)$ or $x$.

Readers can find in  \cite{appell} information about the wide class of standard or Shragin functions, which are not necessarily continuous but enjoy the property of superpositional (or compositional) measurability. More concretely for the purposes in this paper, readers might find useful \cite[proposition 3.2]{cp} to check the validity of $(H_3)-(a)$ with discontinuous $f$ and $\tau$ as we did in example (\ref{ex1}). These technicalities notwithstanding, concrete examples of functions $f$ and $\tau$ which do not satisfy $(H_3)-(a)$ necessarily lean on the axiom of choice in a way one never encounters in applications.
\bigbreak

\noindent
$(H_3)-(b)$ Lots of discontinuities with respect to $x(t)$ are allowed, at least on sufficiently well--behaved sets in the $(t,x)$ plane.

A very interesting feature of $(H_3)-(b)$ is that it allows {\it downwards} jump discontinuities with respect to $x(t)$, something explicitly avoided just a few years ago in the literature on existence of Carath\'eodory solutions for discontinuous first--order differential equations, see \cite{bb}, \cite{bilessh}, \cite{hr}, \cite{hela}. In particular, $(H_3)-(b)$ does not imply continuity with respect to $x(t)$ when combined with $(H_6)$, so $f$ and $\tau$ can be discontinuous with respect to $t$ and to $x(t)$ in Theorems \ref{main2} and \ref{main4}.

\bigbreak

Conditions $(H_3)$ are only used to ensure that the initial value problems (\ref{ivp}) fall inside the scope of Lemma \ref{lemdis}, thus guaranteeing the existence of the extremal solutions satisfying (\ref{greatest}) and (\ref{least}). We could have followed the presentation in \cite{cid}, removing $(H_3)$ and assuming directly that the initial value problems (\ref{ivp}) have extremal solutions fulfilling (\ref{greatest}) and (\ref{least}). In this way, we would leave at the reader's choice which existence result to use in order to get the desired conclusions about existence of solutions to (\ref{ivp}). In other words, our results hold valid if we use another existence result instead of Lemma \ref{lemdis} and if we modify $(H_3)$ accordingly. Some examples of such existence results for initial value problems, which also allow downwards jump discontinuities, are \cite[Theorem 3.1]{bp}, \cite[Theorem 2.1.4]{carlheik}, \cite[Theorem 3.1]{cp2} or \cite[Theorem 3.6]{cp}.

\bigbreak

\noindent
$(H_4)$ Some monotonicity is crucial in our arguments, because they are based on abstract results for monotone operators. This is a price we pay for removing continuity from the list of assumptions. However, thanks to the consideration of the nondecreasing functional term in the differential equation in (\ref{p1}), our results are not so limited by monotonicity as they could seem at first sight. Indeed, consider an equation
\begin{equation}
\label{bv}
x'(t)=g(t,x(t),x(\tau_{t,x(t)})) \quad \mbox{for a.a. $t \in I_+$,}
\end{equation}
where $g:I_+\times \mathbb R^2 \to \mathbb R$ is such that for a.a. $t \in I_+$ and all $x \in \mathbb R$ the mapping $g(t,x,\cdot)$
has bounded variation on compact subsets (for example, if $g(t,x,\cdot)$ is locally Lipschitz continuous). Then there exist functions $g_i:I_+\times \mathbb R^2 \to \mathbb R$ ($i=1,2$) such that $g_1(t,x,\cdot)$ is nondecreasing, $g_2(t,x,\cdot)$ is nonincreasing, and
$g(t,x,\cdot)=g_1(t,x,\cdot)+g_2(t,x,\cdot)$. Therefore equation (\ref{bv}) can be rewritten as
\begin{equation}
\label{bvr}
x'(t)=f(t,x(t),x(\tau_{t,x(t)}),x) \quad \mbox{for a.a. $t \in I_+$,}
\end{equation}
where $f(t,x,y,\gamma)=g_1(t,x,\gamma(\tau_{t,x}))+g_2(t,x,y)$ satisfies the monotonicity properties required in $(H_4)-(a)$, thus falling inside the scope of our theorems as long as monotonicity is concerned. Notice that the concepts of lower and upper solutions depend on the decomposition $g=g_1+g_2$; specifically, $\alpha$ and $\beta$ are lower and upper solutions provided that for a.a. $t \in I_+$ we have
$$\alpha'(t) \le g_1(t,\alpha(t),\alpha(\tau_{t,\alpha(t)}))+g_2(t,\alpha(t),\beta(\tau_{t,\alpha(t)})),$$
and
$$\beta'(t) \ge g_1(t,\beta(t),\beta(\tau_{t,\beta(t)}))+g_2(t,\beta(t),\alpha(\tau_{t,\beta(t)})).$$

\vspace*{0.5cm}
Reformulations of the previous type are also useful for other purposes. Next we show how to avoid checking $(H_3)-(b)$ on some types of discontinuities, namely, those due to nondecreasing terms:

Consider equation (\ref{bv}) again and assume that for a.a. $t \in I_+$ and all $y \in \mathbb R$ the mapping $g(t,\cdot,y)$
has bounded variation on compact subsets. Then $g(t,\cdot,y)=g_1(t,\cdot,y)+g_2(t,\cdot,y)$, where $g_1(t,\cdot,y)$ is nondecreasing and $g_2(t,\cdot,y)$ is nonincreasing. Now equation (\ref{bv}) can be rewritten as (\ref{bvr}) for $f(t,x,y,\gamma)=g_1(t,\gamma(t),y)+g_2(t,x,y)$, and hence discontinuities due to $g_1(t,\cdot,y)$ can occur on arbitrary sets, in the sense that there is no need to check $(H_3)-(b)$ on those discontinuity sets when applying Theorem \ref{main1} or \ref{main3}.

Notice also that $f$ is nonincreasing with respect to $x$, even though $g$ is not, which, to some extent, may simplify finding lower and upper solutions because the results in Section 4 apply, and makes condition $(H_6)$ easier to check as we will show in the corresponding paragraph.

\bigbreak

Lemma 4.2 in \cite{pouso} guarantees that condition (\ref{bv-}) is equivalent to assuming 
\begin{enumerate}
\item function $k$ has bounded variation on $[t_0-\hat r,t_0]$,
\item $k=k_a+k_s$, where $k_a \in AC([t_0-\hat r,t_0])$, $k_s$ is nonincreasing on $[t_0-\hat r,t_0]$, and $k_s'=0$ a.e.,
\item $k_a'(t) \le \hat \psi(t)$ for a.a. $t \in [t_0-\hat r,t_0]$.
\end{enumerate}

Obviously, if $k$ is absolutely continuous on $[t_0-\hat r,t_0]$ then (\ref{bv-}) holds with ``$\le$" replaced by ``$=$" and $\hat \psi$ replaced by $k'$.

\bigbreak
\noindent
$(H_5)$ The one--sided Lipschitz conditions introduced in ($H_5$) become usual Lipschitz conditions when combined with assumptions ($H_4$). Therefore, in the conditions of Theorems \ref{main2} and \ref{main4}, and roughly speaking, $\Lambda$ is a contraction and for a.a. $t \in I_+$, all $x \in [\alpha(t),\beta(t)]$, and all $\gamma \in [\alpha,\beta]^+$ the functions
$$y \in [\alpha(\tau_{t,x,\beta}),\beta(\tau_{t,x,\alpha})] \longmapsto f(t,x, y,\gamma),$$
$$\gamma \in [\alpha,\beta]^+ \longmapsto f(t,x,y,\gamma),$$
and
$$\gamma \in [\alpha,\beta]^+ \longmapsto \int_{t_0}^{\tau(t,x,\gamma)}{\psi_M(s)\,ds}$$
are Lipschitz continuous.

Notice also that (\ref{lipbet}) is fulfilled if $\psi_M$ is essentially bounded and $\tau(t,x,\cdot)$ is Lipschitz continuous on $[\alpha,\beta]^+$.

\bigbreak

\noindent
$(H_6)$ is satisfied, in particular, if the mapping
$$x \in [\alpha(t),\beta(t)] \longmapsto f(t,x,\gamma_1(\tau(t,x,\gamma_2)),\gamma)
$$is nonincreasing. In turn, this is satisfied if (and only if when dropping the $x(t)$ argument in $f$) $f(t,x(t),y,\gamma)$ is nonincreasing in $y$ (already required in $(H_4)-(a)$) and in $x(t)$, and $\tau(t,x(t),x)$ is nondecreasing in $x(t)$. Therefore condition $(H_6)$ can be satisfied even if $f$ or $\tau$ are discontinuous w.r.t. $x(t)$.

\medbreak

$(H_6)$ is also implied by $(H_5)-(a)$ and the following two conditions: for a.a. $t \in I_+$ and all $\gamma_i \in [\alpha,\beta]^+$ such that $\psi_m(t) \le \gamma_i'(t)\le \psi_M(t)$ a.e. on $I_+$ ($\psi_m$ and $\psi_M$ as in $(H_2)$) and $\gamma_i=c_i+k$ on $I_-$ for some constants $c_i \in \mathbb R$ ($i=1,2$) the composition
$$x \in [\alpha(t),\beta] \longmapsto \gamma_1(\tau(t,x,\gamma_2))$$
is Lipschitz continuous (this is valid, in particular, if $\psi_M$ is constant, $k$ is Lipschitz continuous on $I_-$, and $\tau(t,\cdot,\gamma_2)$ is Lipschitz continuous on $[\alpha(t),\beta(t)]$), and there exists $L \in L^1(I_+)$ such that for a.a. $t \in I_+$, all $u \in [\min_{t \in I_{\pm}}\alpha(t),\max_{t \in I_{\pm}}\beta(t)]$ and all $\gamma \in [\alpha,\beta]^+$ we have
$$f(t,x,u,\gamma)-f(t,y,u,\gamma) \le L(t)(y-x)\quad \mbox{whenever $\alpha(t) \le x \le y \le \beta(t)$.}$$
 This shows that neither $f$ nor $\tau$ need be monotone w.r.t. $x(t)$ so that $(H_6)$ is fulfilled.

\begin{remark} The results in this paper can be generalized in the obvious way to the case in which there are several deviating arguments.
\end{remark}

\section{Construction of lower and upper solutions}
In this section we provide two types of sufficient conditions on the data functions for the existence of lower and upper solutions: first, we consider bounded nonlinearities and, second, we study the existence of linear lower and upper solutions for possibly unbounded right--hand sides.

\subsection{Lower and upper solutions for bounded problems}

We start this section with the construction of lower and upper solutions for problem (\ref{p1}) in the case that nonlinearities are bounded on suitable sets. Its proof is easy and we omit it.

\begin{proposition}\label{nosingular} Assume that $f:I_+\times \mathbb R^2 \times {\cal C}(I_{\pm})\to \mathbb R$ is nonincreasing with respect to its third variable and nondecreasing w.r.t. its fourth variable. Assume also that there exist $\lambda_1, \, \lambda_2 \in \mathbb R$ such that $\lambda_1 \le \Lambda(\gamma) \le \lambda_2$ for all $\gamma \in {\cal C}(I_{\pm})$ and let
    $$c_1= \lambda_1+\min_{s \in I_-}k(s),$$
    and
    $$c_2=\lambda_2 + \max_{s \in I_-}k(s).$$
If there exists $\psi \in L^1(I_+)$ such that for a.a. $t \in I_+$, all $x \ge c_2,$ and all $\gamma \in {\cal C}(I_{\pm})$ with $\gamma \ge c_2$ on $I_{\pm}$ we have $$ f(t,x,c_1,\gamma) \leq \psi(t),$$
and for a.a. $t \in I_+$ we have $$0 \leq f(t,c_1,c_2+\|\psi\|_1,c_1),$$
then the functions
\begin{equation*}
\alpha(t)=c_1, \quad t \in I_{\pm},
\end{equation*}
and
\begin{equation*}
\beta(t)=
\left\{
\begin{array}{ll}
c_2, \quad t \in I_-, \\
\\
c_2 + \int_{t_0}^{t} \psi(s) \, ds, \quad t \in I_+,
\end{array}
\right.
\end{equation*}
are, respectively, a lower and an upper solution for problem (\ref{p1}) according to both definitions \ref{subsol} and \ref{subsol2}, and $\alpha \le \beta$ on $I_{\pm}$.
\end{proposition}

\begin{remark} Notice that in the previous proposition no condition about the deviating function $\tau$ is required.
\end{remark}

\subsection{Linear lower and upper solutions}
This section is devoted to highlighting the following facts:
\medbreak
\begin{enumerate}
\item Although finding a pair of lower and upper solutions is not an easy task in general, and it often involves a great deal of mathematical skill, it is also true, as we shall show, that a large class of interesting problems can be studied just with {\it linear} lower and upper solutions.
\medbreak
\item The results in this paper apply for {\it singular problems}, in the sense that the right--hand side in (\ref{p1}) might blow--up at some values of $x$. This is an interesting consequence of the fact that our assumption only have to be satisfied between the given lower and upper solutions.
\end{enumerate}

\medbreak

In this section we consider problem (\ref{p1}) and we assume that the following set of conditions is satisfied:
\medbreak
 \begin{enumerate}
 \item[$(A)$] $0<r \le L$, $k:I_- \to \mathbb R$ is continuous and nondecreasing.
\medbreak
 For the function $\tau:I_+ \times [k(t_0),\infty) \times AC(I_+) \to I_{\pm}$ there exists $\overline{r} \in (0,r]$ such that for all $x \in [k(t_0),\infty)$, all $\gamma \in AC(I_+)$, and a.a. $t \in I_+$ we have
\begin{equation}
\label{re}
\delta(t)=\min \{t_0-\overline{r},t-r\} \le  \tau_{t,x,\gamma}\le t_0+L.  \end{equation}

Let $\Omega=\{ \gamma \in {\cal C}(I_{\pm}) \, : \, \mbox{$\gamma\ge k(t_0-r)$ on $I_{\pm}$}\},$ and let $$f:Dom(f) \subset I_+\times [k(t_0),\infty) \times (k(t_0-r),\infty) \times \Omega \longrightarrow [0, \infty)$$ be non\-in\-crea\-sing w.r.t. its second and third variables and nondecreasing w.r.t. its fourth variable. Moreover, we assume either
$$Dom(f)=I_+\times [k(t_0),\infty) \times (k(t_0-r),\infty) \times \Omega$$
or 
$$Dom(f)=I_+\times (k(t_0),\infty) \times (k(t_0-r),\infty) \times \Omega.$$
\medbreak
The operator $\Lambda:{\cal C}(I_{\pm})\longrightarrow [0,\infty)$ is nondecreasing.

\end{enumerate}

\medbreak

 The first inequality in (\ref{re}) prevents the deviating argument $\tau_{t,x,\gamma}$ from assuming the value $t_0-r$ if $t \neq t_0$, thus allowing solutions to assume singular values at $t_0-r$. This will be made clearer in an example.
 
 \medbreak

For the convenience of the reader we note that the piecewise definition of the function $\delta$ defined in (\ref{re}) is
$$\delta(t) =\left\{ \begin{array}{cl}
  t-r, & \mbox{for $t \in [t_0,t_0+r-\bar r],$} \\
 \\
  t_0-\bar r,  & \mbox{for $t \in [t_0+r-\bar r, t_0+L],$}
  \end{array} \right.$$
  and $\delta(t) \in I_-$ for all $t \in I_+$.
  
  \medbreak

 Next we include the main result in this subsection, which establishes sufficient conditions for the following pair of functions to be, respectively, a lower and an upper solution for (\ref{p1}):
  \begin{equation}
  \label{linalf}
 \alpha(t)= \left\{
 \begin{array}{cl}  k(t), & \mbox{for $t \in I_-,$} \\
\\
m_{\alpha}(t-t_0)+  k(t_0), & \mbox{for $t \in I_+$;} \end{array}
\right.
\end{equation}
\begin{equation}
\label{linbet}
 \beta(t)=\left\{
 \begin{array}{cl}  k(t)+n_{\beta}, & \mbox{for $t \in I_-,$}\\
 \\
 m_{\beta}(t-t_0)+  k(t_0)+n_{\beta}, & \mbox{for $t \in I_+$.}
 \end{array}
 \right.
 \end{equation}

 \begin{proposition}
 \label{proex1}
Assume that conditions $(A)$ hold, and there exist $m_{\alpha}, m_{\beta}, n_{\beta} \in [0, \infty)$ such that $m_{\alpha}\le m_{\beta}$ and that the following conditions are fulfilled:
\begin{eqnarray}
\label{c4}
\Lambda(m_{\beta}L+  k(t_0)+n_{\beta})\le n_{\beta};\\
\label{c3}
f(t,  k(t_0)+n_{\beta}, k(\delta(t)),m_{\beta}L+  k(t_0)+n_{\beta}) \le m_{\beta} \quad \mbox{for a.a. $t \in I_+$ };\\
\label{c2}
m_{\alpha} \le f(t,m_{\alpha}L+k(t_0),m_{\beta}L+  k(t_0)+n_{\beta},k(t_0-r)) \quad \mbox{for a.a. $t \in I_+$}.
\end{eqnarray}

Then the functions $\alpha$ and $\beta$ defined in (\ref{linalf}) and (\ref{linbet}) are, respectively, a lower and an upper solution for (\ref{p1}), according to both Definitions \ref{subsol} and \ref{subsol2}.

In particular, if conditions $(H_2)$, $(H_3)$, and $(H_4)$ in Theorem \ref{main1} hold then (\ref{p1}) has the extremal quasisolutions in $[\alpha,\beta]^+$, and if the remaining conditions in Theorem \ref{main2} hold, then (\ref{p1}) has a unique solution in $[\alpha,\beta]^+$ (and the same is true if the conditions in Theorem \ref{main4} are fulfilled).
\end{proposition}

\noindent {\bf Proof.}  We will only consider Definition \ref{subsol}, as the proof is analogous with Definition \ref{subsol2}. Let us show that $\beta$ is an upper solution (when coupled with $\alpha$). First, for all $t \in I_+$ we have
$$\alpha(\tau_{t,\beta(t),\beta}) \ge \alpha(\delta(t))=
 k(\delta(t)),$$
and thus for a.a. $t \in I_+$ we have
$$f(t,\beta(t),\alpha(\tau_{t,\beta(t),\beta}),\beta) \le f(t,k(t_0)+n_{\beta},k(\delta(t)),m_{\beta}L+  k(t_0)+n_{\beta}) \le m_{\beta}=\beta'(t),$$
by virtue of (\ref{c3}).

Second, inequality (\ref{c4}) guarantees for all $s \in I_-$ that
$$
\Lambda(\beta)+k(s) \le \Lambda(m_{\beta}L+ k(t_0)+n_{\beta})+k(s)\le n_{\beta}+k(s)= \beta (s).
$$

Similar arguments involving (\ref{c2}) show that $\alpha$ is a lower solution.
\qed

Next propositions provide us with a simpler way for checking that linear lower and upper solutions exist. Their proofs are easy so we omit them.

\begin{proposition}
\label{copro1}  (Upper solutions)
Assume that $(A)$ holds. If
\begin{equation}
\label{sob1}
\lim_{z \to +\infty}\frac{f(t,k(t_0)+z,k(\delta(t)),z\,L+  k(t_0)+z)}{z}<1\quad \mbox{almost uniformly in $t \in I_+$}
\end{equation}
and
\begin{equation}
\label{sob2}
\lim_{z \to +\infty}\frac{\Lambda(zL+k(t_0)+z)}{z}<1,
\end{equation}
then (\ref{c4}) and (\ref{c3}) are fulfilled with some $m_{\beta}>0$ and $n_{\beta}=m_{\beta}$.
\end{proposition}

\begin{proposition}
\label{copro2} (Lower solutions)
Assume that $(A)$ holds and that (\ref{c4}) and (\ref{c3}) are satisfied with some $m_{\beta}>0$ and $n_{\beta} \ge 0$.

Then (\ref{c2}) is satisfied with $m_{\alpha}=0$ if $Dom(f)=I_+ \times [k(t_0),\infty) \times (k(t_0)-r,\infty)$, or with some $m_{\alpha} \in (0,m_{\beta})$ if
\begin{equation}
\label{sub}
\lim_{z \to 0^+}\frac{f(t,zL+k(t_0),m_{\beta}L +k(t_0)+n_{\beta},k(t_0-r))}{z}>1 \quad \mbox{almost uniformly in $I_+$.}
\end{equation}
\end{proposition}

\subsubsection{An example}

Let $\{q_n\}_{n \in {\scriptsize \mathbb N}}$ be an enumeration of all rational numbers in $[0,\infty)$ and consider the function
 $$\phi(x)=1-\sum_{\{n \in {\scriptsize \mathbb N} \,:\, q_n \le x\}} 2^{-n} \quad (x \in [0,\infty)),$$
  which is decreasing, continuous at every irrational number in $[0,\infty)$, discontinuous at every rational number in $[0,\infty)$, and $0 <\phi < 1$ on $[0, \infty)$.

  We consider the following general problem with deviated arguments depending functionally on the unknown, possible singularities for $x=0$ and $x=1/2$, and discontinuities with respect to all of its arguments (square brackets $[\cdot]$ mean integer part):
\begin{eqnarray}
\nonumber
x'(t)&=&\varepsilon+f_1(t)+f_2(t)\phi(t+x(t))+f_3(t)\left[\frac{1}{(x(t)-1/2)^{\mu}} \right] \\
&\mbox{\,}&\label{ej1e1}
+\frac{f_4(t)}{x^{\nu}\left(t+\sigma_1(x(t))+\sigma_2\left(x(0)+x(t)+x(1) \right)\right)} \\
\nonumber\\
 \nonumber&\mbox{\,}&+f_5(t)g\left(x\left(t+\sigma_1(x(t))+\sigma_2\left(x(0)+x(t)+x(1) \right)\right)\right), \,\, t \in I_+=[0,1],
 \\
 \nonumber
 \\
\label{ej1e2}
x(t)&=&\varphi \left(\int_0^1{\omega(s) \, x(s)ds}\right)+t+1/2 \quad \mbox{for all $t \in I_-=[-1/2,0]$.}
\end{eqnarray}

Next we introduce some conditions which imply the existence of a couple of linear lower and upper solutions for problem (\ref{ej1e1})--(\ref{ej1e2}).
\begin{proposition}
\label{proloup}
Problem (\ref{ej1e1})--(\ref{ej1e2}) has a couple of lower and upper solutions given by (\ref{linalf}) and (\ref{linbet}) with $m_{\beta}=n_{\beta}$ provided that the following conditions hold:  \begin{enumerate}
 \item[(C1)] (Sufficient for lower and upper solutions) Let $\varepsilon, \mu, \nu \in (0,\infty)$, $f_i \in L^{\infty}(I_+,[0,\infty))$ ($i=1,\dots,5$), $\sigma_i:[0,\infty) \to [-1/4,0]$ ($i=1,2$), $g:[0,\infty) \to [0,\infty)$ is nonincreasing, $\omega$ is nonnegative and integrable on $I_+$, $\varphi:\mathbb R \to [0,\infty)$ is nondecreasing and
\begin{equation}
\label{varfi}
\mbox{either $\omega=0$ a.e. on $I_+$ or $\displaystyle\lim_{z \to +\infty}\frac{\varphi(z)}{z}<\displaystyle\frac{1}{3\int_{I_+}\omega}.$}
\end{equation}
\end{enumerate}
\end{proposition}

\noindent {\bf Proof.}  Problem (\ref{ej1e1})--(\ref{ej1e2}) is a particular case of (\ref{p1}). To see it, we consider $r=1/2$ and $k(t)=t+1/2$ for all $t \in I_-$, and we define $\tau$, $f$ and $\Lambda$ as follows: for all $(t,x,\gamma) \in I_+ \times [0,\infty) \times AC(I_+)$ we define
\begin{equation}
\label{tau}
\tau(t,x,\gamma)=t+\sigma_1(x)+\sigma_2\left(\gamma(0)+\gamma(t)+\gamma(1) \right),
\end{equation}
for $(t,x,y) \in I_+\times (1/2,\infty)\times (0,\infty)$ we define
\begin{eqnarray*}
f(t,x,y)=\varepsilon+f_1(t)+f_2(t)\phi(t+x)+f_3(t)\left[\frac{1}{(x-1/2)^{\mu}} \right]+
\frac{f_4(t)}{y} +f_5(t)g\left(y\right),
 \end{eqnarray*}
 and, finally, $\Lambda(\gamma)=\varphi \left(\int_{I_+}{\omega(s) \gamma(s)ds} \right)$ for all $\gamma \in {\cal C}(I_{\pm})$.

 Notice that $\tau(t,x,\gamma) \in [t-1/2,t]$ for all $(t,x,\gamma) \in I_+ \times [0,\infty) \times AC(I_+)$, $f(t,x,y) \ge \varepsilon$ for all $(t,x,y) \in I_+\times (1/2,\infty)\times (0,\infty)$, and $\Lambda(\gamma) \ge 0$ for all $\gamma \in {\cal C}(I_{\pm})$.

\bigbreak

The conditions on $\sigma_i$ ($i=1,2$) guarantee that (\ref{re}) is satisfied for, say, $\bar r=1/4$, thus $k(\delta(t))=\min \{ t, 1/4 \}$ for $t \in I_+$. The remaining conditions in $(A)$ are immediate consequences of the assumptions in $(C1)$.

\bigbreak

Next we show that the conditions of proposition \ref{copro1} are fulfilled. First, note that there exists $c>0$ such that for a.a. $t \in I_+$ and all $z>0$ we have
$$\frac{f(t,z+1/2,k(\delta(t)))}{z}\le \frac{c}{z}+c\frac{1}{z}\left[\frac{1}{z^{\mu}} \right] +
\frac{c}{z k(\delta(t))} +\frac{c \, g\left(k(\delta(t))\right)}{z},$$
so there exists $d>0$ such that for a.a. $t \in I_+$ and all $z>0$ we have
$$0 \le \frac{f(t,z+1/2,k(\delta(t)))}{z} \le \frac{d}{z}+d\frac{1}{z}\left[\frac{1}{z^{\mu}}\right],$$
which implies (\ref{sob1}). In turn, for $z>1/2$ and $\omega>0$ on a positive measure subset of $I_+$ we have
 $$0 \le \frac{\Lambda(2z+1/2)}{z}\le  \frac{\Lambda(3z)}{z}=\frac{\varphi(3z\int_{I_+}\omega)}{3z\int_{I_+}\omega}3 \int_{I_+}\omega,$$
 so (\ref{sob2}) follows from (\ref{varfi}). Therefore, proposition \ref{copro1} guarantees that there exists $m_{\beta}>0$ such that (\ref{c4}) and (\ref{c3}) are satisfied with $n_{\beta}=m_{\beta}$. Now it suffices to fix $m_{\alpha} \in (0,\varepsilon)$ so that (\ref{c2}) is fulfilled, and then proposition \ref{proex1} applies. (Notice that we have not made any assumption about continuity or monotonicity of the functions $\sigma_i$ yet.)
 \qed

Next we are going to apply Theorem \ref{main1} in order to prove the existence of the extremal quasisolutions to our problem between the lower and upper solutions whose existence we have just proven. To do so, we choose $\hat r=1/2$ to define $[\alpha,\beta]^+$ and we impose the following few additional conditions:
\begin{proposition}
\label{proquasi}
Assume $(C1)$ and
\begin{enumerate}
\item[(C2)] The function
$$t \in I_+ \longmapsto (m_{\alpha} \, t)^{-\mu}f_3(t)+\max\{t^{-\nu},2^{\nu}\}f_4(t)$$
is integrable on $I_+$, $\sigma_1$ is nondecreasing, and $\sigma_2$ is nonincreasing.
\end{enumerate}
Then problem (\ref{ej1e1})--(\ref{ej1e2}) has the extremal quasisolutions in $[\alpha,\beta]^+$.
\end{proposition}

\noindent
\noindent {\bf Proof.}  {\it Claim 1-- $(H_2)$ is fulfilled.} For $t \in (0,1]$ and $x\ge \alpha(t)=m_{\alpha}t+1/2$, we have
\begin{equation}
\label{A}
\left[ \frac{1}{(x-1/2)^{\mu}}\right] \le \frac{1}{(x-1/2)^{\mu}} \le m_{\alpha}^{-\mu}t^{-\mu},
\end{equation}
and for all $\gamma_i \in [\alpha,\beta]^+$ $(i=1,2)$, we have
$$\gamma_1(\tau(t,x,\gamma_2)) \ge \gamma_1(t-1/2)\ge \alpha(t-1/2) \ge \min \{t,1/2\},$$
hence
\begin{equation}
\label{B}
(\gamma_1(\tau(t,x,\gamma_2)))^{-\nu} \le \max\{t^{-\nu},2^{\nu}\}.
\end{equation}

We deduce from (\ref{A}), (\ref{B}), and the properties of the functions $\phi$ and $g$, that for a.a. $t \in I_+$, for all $x \in [\alpha(t),\beta(t)]$, and for all $\gamma_i \in [\alpha,\beta]^+$ $(i=1,2)$, we have
$$0 \le \varepsilon \le f(t,x,\gamma_1(\tau(t,x,\gamma_2)) \le \psi(t),$$
where for a.a. $t \in I_+$
\begin{equation}
\label{psiej}
\psi(t)=\varepsilon+f_1(t)+f_2(t)+m_{\alpha}^{-\mu}t^{-\mu}f_3(t)+\max\{t^{-\nu},2^{\nu}\}f_4(t)+g(0)f_5(t).
\end{equation}
Summing up, condition $(H_2)$ is satisfied as a consequence of the assumptions in $(C2)$.
\medbreak
\noindent
{\it Claim 2 -- $(H_3)-(a)$ is fulfilled.} We have to prove that for all $x \in [\alpha(0),\beta(1)]$, and for all $\gamma_i \in [\alpha,\beta]^+$ $(i=1,2)$ the composition
$$t \in \{t \in I_+ \, : \, \alpha(t) \le x \le \beta(t)\} \longmapsto f(t,x,\gamma_1(\tau(t,x,\gamma_2))),$$
is measurable. This is trivial for $x=\alpha(0)$ or $x=\beta(1)$, because the corresponding domains are null--measure sets, so we restrict our attention to $x \in (\alpha(0),\beta(1))$. Remember that monotone functions are Borel--measurable, and that sums or compositions of Borel--measurable functions are Borel--measurable. For $t \in (0,1]$ we have
 $$\gamma_1(\tau(t,x,\gamma_2)) \ge \gamma_1(t-1/2)>0,$$
and therefore the mappings
 $$t   \longmapsto \frac{1}{{\gamma_1}^{\nu}\left(t+\sigma_1(x)+\sigma_2\left(\gamma_2(0)+\gamma_2(t)+\gamma_2(1) \right)\right)}$$
 and
 $$t   \longmapsto g\left(\gamma_1\left(t+\sigma_1(x)+\sigma_2\left(\gamma_2(0)+\gamma_2(t)+\gamma_2(1) \right)\right)\right)$$
 are Borel--measurable, which guarantees that
\begin{eqnarray*}
\nonumber
f(t,x,\gamma_1(\tau_{t,x,\gamma_2}))&=&\varepsilon+f_1(t)+f_2(t)\phi(t+x)+f_3(t)\left[\frac{1}{(x-1/2)^{\mu}} \right] \\
&\mbox{\quad}& +\frac{f_4(t)}{ {\gamma_1}^{\nu}\left(\tau_{t,x,\gamma_2} \right)  } +f_5(t)g\left( \gamma_1\left( \tau_{t,x,\gamma_2} \right)    \right),
 \end{eqnarray*}
is measurable on the domain $\{t \in I_+ \, : \, \alpha(t) \le x \le \beta(t)\}$.

\medbreak

\noindent
{\it Claim 3 -- $(H_3)-(b)$ is fulfilled.} For a.a. $t \in I_+$ and all $\gamma_i \in [\alpha,\beta]^+$ ($i=1,2$) the mapping
$$x \in [\alpha(t),\beta(t)] \longmapsto f(t,x,\gamma_1(\tau(t,x,\gamma_2))),$$
is continuous on $[\alpha(t),\beta(t)]$ except, at most, for
\begin{enumerate}
\item $x=q_n-t$ for some $n \in \mathbb N$;
\item $x=n^{-1/\mu}+1/2$ for some $n \in \mathbb N$; or
\item $x=v_n$ for some $n \in \mathbb N$, $\{v_n\}_{n \in {\scriptsize \mathbb N}}$ being the sequence of all discontinuity points of $\sigma_1$; or
    \item $x=w_n$ for some $n \in \mathbb N$, $\{w_n\}_{n \in {\scriptsize \mathbb N}}$ being the sequence of all discontinuity points of the nonincreasing function
        $$x \in [0,\infty) \longmapsto g(\gamma_1(t+\sigma_1(x)+\sigma_2(\gamma_2(0)+\gamma_2(t)+\gamma_2(1)))).$$
\end{enumerate}
Therefore all possible discontinuity points of $f(t,\cdot,\gamma_1(\tau(t,\cdot,\gamma_2)))$ lie over countable many lines on the $(t,x)$ plane which have nonpositive slope. On the other hand we know that $f \ge \varepsilon>0$ on its domain, and therefore condition $(H_3)-(b)$ is satisfied.

Notice that $(H_4)$ is an obvious consequence of $(C1)$ and $(C2)$, so we are in a position to apply Theorem \ref{main1} to ensure that the problem has the extremal quasisolutions in $[\alpha,\beta]^+$.
\qed

\bigbreak

Finally, we are going to introduce a new set of conditions so that our problem fulfills all the assumptions in Theorem \ref{main3}, thus assuring the existence of a unique solution in $[\alpha,\beta]^+$.

\begin{proposition}
\label{prouni}
We assume $(C1)$, $(C2)$, and
\begin{enumerate}
\item[(C3)] The function
$$t \in I_+ \longmapsto t^{-\nu-1}f_4(t)$$
is integrable on $I_+$, there exists $L_g \ge 0$ such that for $0 \le u\le v$ we have
$$g(u)-g(v) \le L_g (v-u),$$
the function $\psi$ given in (\ref{psiej}) is essentially bounded, there exists $L_{\sigma_2} \ge 0$ such that for $0 \le x \le y$ we have
$$\sigma_2(x)-\sigma_2(y) \le L_{\sigma_2}(y-x),$$
and there exists $L_{\varphi} \ge 0$ such that for $0 \le x \le y$ we have
$$\varphi(y)-\varphi(x) \le L_{\varphi} (y-x).$$
\end{enumerate}
Then problem (\ref{ej1e1})--(\ref{ej1e2}) has a unique solution in $[\alpha,\beta]^+$ provided that
\begin{eqnarray}
\label{cond}
L_{\varphi} \int_0^1{\omega(s)ds}
+ (1+ 3\|\tilde\psi\|_{\infty}L_{\sigma_2})\int_{0}^{1}{\left(\nu \max\{s^{-1-\nu}2^{1+\nu}\}f_4(s)+L_gf_5(s)\right)ds} <1,
\end{eqnarray}
where $\tilde \psi=1$ on $I_-$ and $\tilde \psi=\psi$ on $I_+$, $\psi$ given in (\ref{psiej}).
\end{proposition}

\noindent
\noindent {\bf Proof.}  {\it Claim 1 -- $(H_5)-(a)$ is fulfilled.} Notice that for all $t \in (0,1]$ and all $x \in [\alpha,\beta]$ we have $0<\alpha(t-1/2) \le \alpha(\tau(t,x,\beta))$, thus by the mean value Theorem for $\alpha(\tau(t,x,\beta)) \le u \le v \le \beta(\tau(t,x,\alpha))$ there is some $z \in [u,v]$ such that
$$u^{-\nu}-v^{-\nu}=\nu z^{-1-\nu}(v-u) \le \nu \alpha^{-1-\nu}(t-1/2).$$
Moreover, for all $t \in (0,1/2]$ we have $\alpha(t-1/2)=t$, and for all $t \in [1/2,1]$ we have $\alpha(t-1/2)=m_{\alpha}(t-1/2)+1/2 \ge 1/2$, so for all $t \in (0,1]$ we have
$$\alpha^{-1-\nu}(t-1/2) \le \max\{t^{-1-\nu},2^{1+\nu}\}.$$
Summing up, for all $t \in (0,1]$, all $x \in [\alpha(t),\beta(t)]$, and $\alpha(\tau(t,x,\beta)) \le u \le v \le \beta(\tau(t,x,\alpha))$, we have
$$f(t,x,u)-f(t,x,v) \le\{ \nu \max\{ t^{-1-\nu}2^{1+\nu}\}f_4(t)+L_gf_5(t)\}(v-u),$$
and therefore $(H_5)-(a)$ is fulfilled by virtue of $(C3)$ and choosing
$$L_1(t)=  \nu \max\{ t^{-1-\nu}2^{1+\nu}\}f_4(t)+L_gf_5(t) \quad \mbox{for a.a. $t \in I_+$.}$$

\medbreak

\noindent
{\it Claim 2 -- $(H_5)-(b)$ is fulfilled.} Following the statement of Theorem \ref{main3} we can define
$\tilde \psi=1$ on $I_-$, and $\tilde \psi=\psi$ on $I_+$, $\psi$ given in (\ref{psiej}). Condition $(C3)$ guarantees that $\tilde \psi \in L^{\infty}(I_{\pm})$. Now, we have for a.a. $t \in I_+$, all $x\in [\alpha(t),\beta(t)]$, and all $\gamma_i \in [\alpha,\beta]^+$ $(i=1,2)$, $\gamma_1 \le \gamma_2$ on $I_{\pm}$, that
\begin{align*}
\int_{\tau_{t,x,\gamma_2}}^{\tau_{t,x,\gamma_1}}{\tilde \psi(s)ds} & \le \|\tilde \psi\|_{\infty}
(\sigma_2(\gamma_1(0)+\gamma_1(t)+\gamma_1(1))  -\sigma_2(\gamma_2(0)+\gamma_2(t)+\gamma_2(1)))\\
& \le 3\|\tilde \psi\|_{\infty}L_{\sigma_2}\max_{s \in I_{\pm}}(\gamma_2(s)-\gamma_1(s)),
\end{align*}
so $(H_5)-(b)$ holds with $L_3(t)=3\|\tilde \psi\|_{\infty}L_{\sigma_2}$ on $I_+$.

$(H_5)-(c)$ is satisfied with $\lambda$ replaced by $L_{\varphi}\int_{I_+}{\omega}$, and the monotonicity properties of the involved functions yield that $(H_6)$ holds with $L_4=0$. These proofs are easy and we omit it.

Theorem \ref{main3} ensures now that problem (\ref{ej1e1})--(\ref{ej1e2}) has a unique solution in $[\alpha,\beta]^+$ provided that (\ref{cond}) holds.
\qed

We emphasize the fact that $\sigma_i$ $(i=1,2)$ might be discontinuous on dense subsets of $[0,\infty)$, and that the specific form of the functional dependences in our problem, namely,
$$x(0)+x(t)+x(1) \quad \mbox{and} \quad \int_0^1{\omega(s)x(s)ds},$$
is unessential, their only important feature is nondecreasingness with respect to $x(\cdot)$.

The results in this section are very general, as we want to emphasize by sorting out a remarkable particular case of (\ref{ej1e1})--(\ref{ej1e2}), which satisfies the conditions in proposition \ref{prouni}.

\begin{example}
Let $C_i \subset [0,1]$ be Cantor sets with positive Lebesgue measures $M_i$ ($i=1,\dots,5$), and let $\chi_{i}$ stand for their respective characteristic functions. We assume that $C_4 \subset [0,1/2]$ only for simplicity in computations.

The following singular and discontinuous problem with state--dependent deviations, where \\
$\sigma:[0,\infty) \to [-1/2,0]$ is nondecreasing (not necessarily continuous), has a unique solution between certain linear lower and upper solutions provided that $M_4+M_5<1$:

\begin{eqnarray}
\nonumber
x'(t)&=&1+\chi_1(t)+\chi_2(t)\phi(t+x(t))+t^3\chi_3(t)\left[\frac{1}{(x(t)-1/2)^{2}} \right] \\
&\mbox{\,}&\label{ej1e1+}
+\frac{t^5\chi_4(t)}{4x^{4}\left(t+\sigma(x(t))\right)} +\frac{\chi_5(t)}{x\left(t+\sigma(x(t))\right)+1}, \quad t \in I_+=[0,1],
 \\
 \nonumber
 \\
\label{ej1e2+}
x(t)&=&t+1/2 \quad \mbox{for all $t \in I_-=[-1/2,0]$.}
\end{eqnarray}

Notice that the differential equation is singular at $x=1/2$ (due to the fourth term in the right--hand side of (\ref{ej1e1+})), and it is also singular at $x=0$ provided that $\sigma(1/2)=-1/2$ (due to the fifth term). The initial condition (\ref{ej1e2+}) forces the unique solution to assume both singular values at $t=0$ if $\sigma(1/2)=-1/2$.

\end{example}

\section{Sign--changing nonlinearities}
Dropping the functional dependence in the deviating function $\tau$ allows us to consider diffe\-ren\-tial equations with sign--changing nonlinearities, nonmonotone initial data, and nonmonotone lower and upper solutions.

In this section we are concerned with the following kind of problems:
\begin{equation}\label{po11}
\left\{
\begin{array}{ll}
x'(t)=f(t,x(t),x(\tau_{t,x(t)}),x)   \quad \mbox{for a.a. $t \in I_+=[t_0,t_0+L],$} \\
\\
x(t)=\Lambda(x)+k(t) \quad \mbox{for all $t \in I_-=[t_0-r,t_0].$}
\end{array}
\right.
\end{equation}

Notice that the previous two different notions of quasisolutions given in Definitions \ref{quasi} and \ref{quasi2} coincide when adapted to the particular case of problem (\ref{po11}), and the same is true for the concepts of lower and upper solutions established in Definitions \ref{subsol} and \ref{subsol2}.

\medbreak

First we state a new result on existence of extremal quasisolutions for problem (\ref{po11}) and we give a sketch of its proof.

\begin{theorem}\label{mainpo1} Suppose that

\noindent
either $r=0$ and $k$ is a fixed real number, or $r>0$ and $k:I_- \longrightarrow \mathbb R$ is continuous on $I_-$.

Assume

\medbreak

\begin{enumerate}
\item[($\hat H_1$)]{(Lower and upper solutions) There exist $\alpha,\beta$, a lower and an upper solution to problem (\ref{po11}), respectively, such that $\alpha \leq \beta$ on $I_{\pm}$,}

\medbreak

\item[$(\hat H_2)$] ($L^1$ bound) There exists $\psi \in L^1(I_+)$ such that for a.a. $t \in I_+$, all $x \in [\alpha(t),\beta(t)]$, and all $\gamma:I_{\pm} \longrightarrow \mathbb R$ such that $\alpha \leq \gamma \leq \beta$ on $I_{\pm},$ we have $|f(t,x,\gamma (\tau(t,x)),\gamma)|\le \psi(t)$,
    \end{enumerate}

    \medbreak

    \noindent
and let $[\alpha,\beta]=\{\gamma \in {\cal C}(I_{\pm})  \, : \, \mbox{$\alpha \leq \gamma \leq \beta$ on $I_{\pm}$, $\gamma_{|I_+} \in AC(I_+)$,  and $|\gamma'| \leq \psi$ a.e. on $I_+$} \}.$

\medbreak

Suppose also that the following conditions hold:

\medbreak

\begin{enumerate}

\item[($\hat H_3$)] (Discontinuous and non--monotone dependences on $t$ and $x(t)$)
\begin{enumerate}
\item[(a)] (Measurability w.r.t. $t$)

For all $x \in [\min_{t \in I_+}\alpha(t),\max_{t \in I_+}\beta(t)]$ and all $\gamma_i \in [\alpha,\beta] $ ($i=1,2$), the compositions $t \in \{s \in I_+ \, : \, \alpha(s) \leq x \leq \beta(s)\} \longmapsto f(t,x,\gamma_1(\tau(t,x)),\gamma_2)$, and

    \noindent
    $t \in I_+ \longmapsto f(t,\gamma_2(t),\gamma_1(\tau(t,\gamma_2(t))),\gamma_2)$ are measurable;

    \medbreak

    \item[(b)] (Admissible discontinuities w.r.t. $x(t)$)

    For a.a. $t \in I_+$ and all $\gamma_i  \in [\alpha,\beta] $ ($i=1,2$), the function $f(t,\cdot,\gamma_2(\tau(t,\cdot )),\gamma_1)$ is continuous on $[\alpha(t),\beta(t)]
\backslash K(t)$, where $
K(t)= \cup_{n=1}^{\infty}K_n(t)$ and may depend on the choice of $\gamma$,
and for
each $n \in \mathbb N$ and each $x \in K_n(t)$ we have
$$
\bigcap_{\varepsilon >0}
\overline{{\rm co}} f(t, x + \varepsilon B,\gamma_2 (\tau_{t,x+\varepsilon B}),\gamma_1)
\bigcap DK_n(t,x)(1)  \subset
\{f(t,x,\gamma_2 (\tau_{t,x}),\gamma_1)\},
$$
where $B$ and $\overline{{\rm co}}$ are as in Lemma \ref{lemdis}.

    \end{enumerate}

    \medbreak

\item[($\hat H_4$)] (Monotone functional dependences)

\begin{enumerate}

\item[(a)] For a.a. $t \in I_+$, for all $x \in [\alpha(t),\beta(t)]$, and all $\gamma \in [\alpha,\beta]$ the function $f(t,x,\cdot,\gamma)$ is nonincreasing on $[\min_{t \in I_{\pm}}\alpha(t),\max_{t \in I_{\pm}}\beta(t)]$; and for a.a. $t \in I_+$, for all $x \in [\alpha(t),\beta(t)]$, and for all $y \in [\min_{t \in I_{\pm}}\alpha(t),\max_{t \in I_{\pm}}\beta(t)]$ the operator $f(t,x,y,\cdot)$ is nondecreasing on $[\alpha,\beta]$, i.e., for $\gamma_i \in [\alpha,\beta]$ ($i=1,2$) the relation $\gamma_1 \le \gamma_2$ on $I_{\pm}$ implies $f(t,x,y,\gamma_1) \le f(t,x,y,\gamma_2)$;

\medbreak

\item[(b)] $\Lambda$ is nondecreasing on $  [\alpha,\beta]  $.
\end{enumerate}
\end{enumerate}
\medbreak
Then problem (\ref{po11}) has the extremal quasisolutions in $v_*, w^* \in [\alpha,\beta]$, which satisfy
$$
(v_*,w^*)=\min_{\preceq} \{(v,w) \, : \, (w,v) \mbox{ are lower and upper solutions in $[\alpha,\beta]$ for } (\ref{po11})\}.$$

\end{theorem}

\noindent {\bf Proof.}  It suffices to follow step by step the proof of Theorem \ref{main1} with the space  $X=\mathcal{C}(I_{\pm})$, endowed with its usual metric and pointwise ordering, the subset $\hat Y=\{ \gamma \in X \, : \, \mbox{$\gamma_{|I_+} \in AC(I_+), \ |\gamma'| \leq \psi$ a.e.}\},$ the ordered interval  $[\alpha,\beta]$ introduced in the statement, and the multivalued operator ${\cal A}:[\alpha,\beta]  \times [\alpha,\beta]  \longrightarrow 2^{[\alpha,\beta]} \setminus \emptyset$ defined as follows: for each pair $(\gamma_1,\gamma_2) \in [\alpha,\beta] \times [\alpha,\beta]$ let ${\cal A}(\gamma_1,\gamma_2)$ be the set of solutions in $[\alpha,\beta]$ to the initial value problem
\begin{equation}\label{ivpaux2}
z'(t)=f(t,z(t),\gamma_2(\tau_{t,z(t)}),\gamma_1) \mbox{ a.e. on $I_+$}, \,\, z=\Lambda({\gamma_1})+k \, \, \mbox{on $I_-.$}
\end{equation}
The most important difference with respect to the proofs of Theorems \ref{main1} or \ref{main2} is that the mixed monotonicity of operators $A_{\pm}$ does not depend on any monotonicity property of the functions in $[\alpha,\beta]$ (see (\ref{des1})).
\qed

Now we establish a uniqueness result for (\ref{po11}).

\begin{theorem} \label{mainpo2} Suppose that all the assumption in Theorem \ref{mainpo1} hold, and assume that the following group of conditions is satisfied:

\medbreak

\begin{enumerate}

    \item[($\hat H_5$)] (One--sided Lipschitz functional dependences)
    \begin{enumerate}
    \item[(a)] Let $L_1, \, L_2 \in L^1(I_+,[0,\infty))$ be such that for a.a. $t \in I_+$, all $x \in [\alpha(t),\beta(t)]$, and all $\gamma_i \in [\alpha,\beta]$ ($i=1,2$) the relations $\alpha(\tau_{t,x})\le u \le v \le \beta(\tau_{t,x})$ and $\gamma_1 \le \gamma_2$ on $I_{\pm}$ imply
  \begin{equation}\label{hpo4}
f(t,x,u,\gamma_2)-f(t,x,v,\gamma_1) \le L_1(t)(v-u)+L_2(t) \max_{s \in I_{\pm}}(\gamma_2(s)-\gamma_1(s));
\end{equation}

\medbreak

    \item[(b)] There exists $\lambda \in [0,1)$ such that for $\gamma_i \in [\alpha,\beta]$ ($i=1,2$), $\gamma_1 \le \gamma_2$, we have
    $$\Lambda({\gamma_2})-\Lambda({\gamma_1}) \le \lambda
    \max_{s \in I_{\pm}}(\gamma_2(s)-\gamma_1(s)).$$
    \end{enumerate}

    \medbreak

    \item[$(\hat H_6)$] (One--sided Lipschitz condition w.r.t. $x(t)$)

    There exists $L_3 \in L^1(I_+,[0,\infty))$ such that for a.a. $t \in I_+$ and all $\gamma_i \in [\alpha,\beta]$ ($i=1,2$) such that $\gamma_i=c_i+k$ on $I_-$ for some constants $c_i \in \mathbb R$, the relation $\alpha(t) \le x \le y \le \beta(t)$ implies
  \begin{equation}\label{hpo41}
f(t,y,\gamma_2(\tau(t,y)),\gamma_1)-f(t,x,\gamma_2(\tau(t,x)),\gamma_1) \le L_3(t)(y-x).
\end{equation}
\end{enumerate}
Then problem (\ref{po11}) has a unique solution in $[\alpha, \beta]$ provided that  \begin{equation}
\label{small2}
\lambda + \int_{t_0}^{t_0+L}{L_1(s)ds}+\int_{t_0}^{t_0+L}{L_2(s)ds}+ \int_{t_0}^{t_0+L}{L_3(s)ds}<1.
\end{equation}
\end{theorem}
\noindent
\noindent {\bf Proof.}  Repeat the proof of Theorem \ref{main2} and notice that the computations are easier in this case and do not depend on monotonicity properties of $v_*$ or $w^*$. Indeed, for a.a. $t \in I_+$ we have
\begin{align*}
p'(t)&=f(t,w^*(t),v_*(\tau_{t,w^*(t) }),w^*)-f(t,v_*(t),v_*(\tau_{t,v_*(t) }),w^*)
\\
& \mbox{\quad}
+f(t,v_*(t),v_*(\tau_{t,v_*(t) }),w^*)-f(t,v_*(t),w^*(\tau_{t,v_*(t) }),v_*) \\
&\le L_3(t)(w^*(t)-v_*(t)) +L_1(t)(w^*(\tau_{t,v_*(t) })-v_*(\tau_{t,v_*(t) }))+L_2(t) \max_{s \in I_{\pm}}p(s) \\
 & \le (L_1(t)+L_2(t)+L_3(t))\max_{s \in I_{\pm}}p(s).
\end{align*}
\qed

\begin{remark} Note that Theorems \ref{mainpo1} and \ref{mainpo2} improve the results in \cite{jan}, even in the continuous case, because $f$ can change sign and $\tau(t,\cdot)$ can be nonmonotone.
\end{remark}

We end this section with an example that, as far as we know, is not covered by any result in the literature.

\begin{example}\label{ej2} Consider the problem

\begin{equation}\label{ex2}
x'(t)=-f_1(t)x^{n}(\sin(t^{-1} + f_2(t)x)), \  \mbox{for a.a.  } t \in [0,1], \quad x(t)=k(t), \ \mbox{for }  t \in [-1,0],
\end{equation}
where
$$
k(t)= \left\{
\begin{array}{ll}
t^2 \, \sin \left( \dfrac{1}{t} \right), \  &-1 \leq t <0, \\
0, \ &t=0,
\end{array}
\right.
$$

\noindent $f_1$ is a nonnegative integrable function with $\|f_1\|_1 \leq 1$, $f_2$ and $f_1 f_2$ are integrable, and $n \in \mathbb N$. Finally we assume that the mapping $$t \in [0,1] \longrightarrow \left( \int_0^t f_1(s) \, ds \right) ^{n}$$ is lipschitzian. \\
Notice that in this example both the start function and the function de\-fi\-ning the equation have nonconstant sign. \\

\noindent
{\it Case 1-- $n$ odd.} The reader can check that conditions $(\hat{H}_1)-(\hat{H}_6)$ hold with:
\begin{enumerate}
\item  $\alpha(t)=-\int_0^t f_1(s) \, ds =-\beta(t)$, for $t \in [0,1],$ and $\alpha=\beta=k$ on $[-1,0]$;
\item $\psi=f_1$;
\item $L_1=n f_1$, \ $L_2=0$, \ $L_3=\|(\beta^{n})'\|_{L^{\infty}(-1,1)} f_1 |f_2|$ and $\lambda =0$.
\end{enumerate}

Then Theorem \ref{mainpo2} guarantees that problem (\ref{ex2}) has an unique solution in the functional interval $[\alpha,\beta]$ provided that $$n\|f_1\|_1 + \|(\beta^{n})'\|_{L^{\infty}(-1,1)} \|f_1 f_2\|_1 < 1.$$

For the convenience of the reader, we check $(\hat{H}_6)$: for a.a. $t \in I$, all $\gamma \in [\alpha,\beta]$ such that $\gamma=k$ on $I_-$, and $\alpha(t)\leq x \leq y \leq \beta(t)$ we have
$$
f_1(t) [\gamma^{n}(\sin(t^{-1} + f_2(t)x)) - \gamma^{n}(\sin(t^{-1} + f_2(t)y))] = f_1(t) \int_{\sin(t^{-1} + f_2(t)y)}^{\sin(t^{-1} + f_2(t)x)} (\gamma^{n})'(s) \, ds,
$$
where $(\gamma^{n})'=(k^{n})'=(\beta^{n})'$ a.e. on $[-1,0]$ and for a.a. $s \in [0,1]$ we have
$$
(\gamma^{n})'(s)=n (\gamma^{n -1})'(s) \gamma'(s) \leq (\beta^{n})'(s).$$
Therefore, for a.a. $t \in I$, all $\gamma \in [\alpha,\beta]$ such that $\gamma=k$ on $I_-$, and $\alpha(t)\leq x \leq y \leq \beta(t)$ we have
$$
f_1(t) [\gamma^{n}(\sin(t^{-1} + f_2(t)x)) - \gamma^{n}(\sin(t^{-1} + f_2(t)y))]  \leq \|(\beta^{n})'\|_{L^{\infty}(-1,1)} f_1(t) |f_2(t)| (y-x),
$$

\noindent thus proving that $(\hat{H}_6)$ holds with $L_3=\|(\beta^{n})'\|_{L^{\infty}(-1,1)} f_1 |f_2|$.

\medbreak

\noindent
{\it Case 2 -- $n$ even.} In this case the differential equation is not monotone with respect to the functional argument, so we have to rewrite it in a way that our theorems apply. To do so, we follow some ideas that we described in Section 3.4, and we express
$$z \in \mathbb R \longmapsto z^n=g_1(z)+g_2(z),$$
where $g_1(z)=z^n \chi_{(-\infty,0)}(z)$ is nonincreasing and $g_2(z)=z^n \chi_{[0,\infty)}(z)$ is nondecreasing. Now for $(t,x,y,\gamma) \in (0,1] \times {\mathbb R}^2 \times {\cal C}[-1,1]$ we define
$$\tau(t,x)=\sin(t^{-1} + f_2(t)x),$$
and
$$f(t,x,y,\gamma)=-f_1(t) \left(g_1(\gamma(\tau_{t,x}))+g_2(y) \right),$$
so the differential equation in (\ref{ex2}) becomes
$$x'(t)=f(t,x(t),x(\tau_{t,x}),x) \quad \mbox{for a.a. $t \in I_+$.}$$

Let us prove $\alpha(t)=-\int_0^t f_1(s) \, ds =-\beta(t)$, for $t \in [0,1],$ and $\alpha=\beta=k$ on $[-1,0]$, define lower and upper solutions provided that $\|f_1 \|_1 \le 2^{-1/n}$. Indeed, for a.a. $t \in [0,1]$ such that
$\tau_{t,\alpha(t)}>0$ we have
$$
g_1(\alpha(\tau_{t,\alpha(t)}))+g_2(\beta(\tau_{t,\alpha(t)}))=\alpha^n(\tau_{t,\alpha(t)})+\beta^n(\tau_{t,\alpha(t)})  \le 2 \|f_1\|_1^n \le 1,
$$
and if $\tau_{t,\alpha(t)} \le 0$ then $
g_1(\alpha(\tau_{t,\alpha(t)}))+g_2(\beta(\tau_{t,\alpha(t)}))=k^n(\tau_{t,\alpha(t)}) \le 1$. Hence, for a.a. $t \in [0,1]$ we have
$$\alpha'(t)=-f_1(t) \le f(t,\alpha(t),\beta(\tau_{t,\alpha(t)}),\alpha).$$
In an analogous way one can prove that $\beta'(t) \ge f(t,\beta(t),\alpha(\tau_{t,\beta(t)}),\beta)$ for a.a. $t \in [0,1]$.
Then, hypotheses $(\hat{H}_1)$-$(\hat{H}_6)$ in Theorem \ref{mainpo2} hold with
\begin{enumerate}
\item $\psi=f_1$;
\item $L_1=L_2=n f_1$, \ $L_3=2\|(\beta^{n})'\|_{L^{\infty}(-1,1)} f_1 |f_2|$ and $\lambda =0$.
\end{enumerate}
Therefore, problem (\ref{ex2}) has a unique solution in $[\alpha,\beta]$ provided that
$$n\|f_1\|_1 + \|(\beta^{n})'\|_{L^{\infty}(-1,1)} \|f_1 f_2\|_1 < 1/2.$$

\end{example}

\subsection{Extremal solutions for nondecreasing functional equations}

As a consequence of Theorem \ref{mainpo1} we have an existence result for problems of the type
\begin{equation}\label{qo11}
\left\{
\begin{array}{ll}
x'(t)=f(t,x(t),x)   \quad \mbox{for a.a. $t \in I_+=[t_0,t_0+L],$} \\
\\
x(t)=\Lambda(x)+k(t) \quad \mbox{for all $t \in I_-=[t_0-r,t_0],$}
\end{array}
\right.
\end{equation}
where $f$ is nondecreasing with respect to its third argument. In this case lower and upper solutions as defined in Definitions \ref{subsol} or \ref{subsol2} become uncoupled, and the extremal quasisolutions in $Y \subset {\cal C}(I_{\pm})$ turn into the least solution and the greatest one in $Y$. Gathering all this information we have the following corollary.

\begin{corollary}\label{mainpo111} Suppose that

\noindent
either $r=0$ and $k$ is a fixed real number, or $r>0$ and $k:I_- \longrightarrow \mathbb R$ is continuous on $I_-$.

Assume

\medbreak

\begin{enumerate}
\item[($\tilde H_1$)]{(Lower and upper solutions) There exist $\alpha,\beta$, a lower and an upper solution to problem (\ref{qo11}), respectively, such that $\alpha \leq \beta$ on $I_{\pm}$,}

\medbreak

\item[$(\tilde H_2)$] ($L^1$ bound) There exists $\psi \in L^1(I_+)$ such that for a.a. $t \in I_+$, all $x \in [\alpha(t),\beta(t)]$, and all $\gamma:I_{\pm} \longrightarrow \mathbb R$ such that $\alpha \leq \gamma \leq \beta$ on $I_{\pm},$ we have $|f(t,x,\gamma)|\le \psi(t)$,
    \end{enumerate}

    \medbreak

    \noindent
and let $[\alpha,\beta]=\{\gamma \in {\cal C}(I_{\pm})  \, : \, \mbox{$\alpha \leq \gamma \leq \beta$ on $I_{\pm}$, $\gamma_{|I_+} \in  AC(I_+)$, and $|\gamma'| \leq \psi$ a.e. on $I_+$} \}.$

\medbreak

Suppose also that the following conditions hold:

\medbreak

\begin{enumerate}
\item[($\tilde H_3$)] (Discontinuous and non--monotone dependences on $t$ and $x(t)$)
\begin{enumerate}
\item[(a)] (Measurability w.r.t. $t$)

For all $x \in [\min_{t \in I_+}\alpha(t),\max_{t \in I_+}\beta(t)]$ and all $\gamma \in [\alpha,\beta] $, the compositions

\noindent $t \in \{s \in I_+ \, : \, \alpha(s) \leq x \leq \beta(s)\} \longmapsto f(t,x,\gamma)$, and $t \in I_+ \longmapsto f(t,\gamma(t),\gamma)$ are mea\-su\-rable;

    \medbreak

    \item[(b)] (Admissible discontinuities w.r.t. $x(t)$)

    For a.a. $t \in I_+$ and all $\gamma  \in [\alpha,\beta]$, the function $f(t,\cdot,\gamma)$ is continuous on $[\alpha(t),\beta(t)] \backslash K(t)$, where $
K(t)= \cup_{n=1}^{\infty}K_n(t)$ and may depend on the choice of $\gamma$,
and for
each $n \in \mathbb N$ and each $x \in K_n(t)$ we have
$$
\bigcap_{\varepsilon >0}
\overline{{\rm co}} f(t, x + \varepsilon B,\gamma)
\bigcap DK_n(t,x)(1)  \subset
\{f(t,x,\gamma)\},
$$
where $B$ and $\overline{{\rm co}}$ are as in Lemma \ref{lemdis}.

    \end{enumerate}

    \medbreak

\item[($\tilde H_4$)] (Monotone functional dependences)

\begin{enumerate}

\item[(a)] For a.a. $t \in I_+$ and for all $x \in [\alpha(t),\beta(t)]$ the operator $f(t,x,\cdot)$ is nondecreasing on $[\alpha,\beta]$, i.e., for $\gamma_i \in [\alpha,\beta]$ ($i=1,2$) the relation $\gamma_1 \le \gamma_2$ on $I_{\pm}$ implies $f(t,x,\gamma_1) \le f(t,x,\gamma_2)$;

\medbreak

\item[(b)] $\Lambda$ is nondecreasing on $ [\alpha,\beta] $.
\end{enumerate}
\end{enumerate}
Then problem (\ref{qo11}) has the extremal solutions in $[\alpha,\beta] $.
\end{corollary}

\begin{remark}
Corollary \ref{mainpo111} is the unique result in this paper which ensures existence of solutions but not uniqueness. In contrast, no Lipschitz--type assumption is needed, and Corollary \ref{mainpo111} yields an existence result for problems with periodic boundary conditions, which do not fulfill the contractivity conditions (\ref{small}) or (\ref{small2}).
\end{remark}

 \end{document}